\documentclass[11pt]{article}

\usepackage{relsize}
\usepackage{amssymb}
\usepackage{amsthm}
\usepackage{mathtools}

\usepackage{stackengine}

\newcommand{\preceqq}{\mathrel{\raisebox{-3.8pt}{\stackanchor[1pt]{$\prec$}{$=$}}}}

\usepackage{afterpage}
\usepackage{makebox}

\DeclareMathSymbol{\alphld}{\mathalpha}{AMSa}{"45}
\DeclareMathSymbol{\Alpha}{\mathalpha}{operators}{"41}

\usepackage{setspace}
\usepackage[usenames,dvipsnames]{color}
\usepackage{xspace}
\usepackage{enumitem}
\newlength{\temp}

\newcommand{\medcup}{{\textstyle\bigcup}}
\newcommand{\medcap}{{\textstyle\bigcap}}
\renewcommand{\preceq}{\preccurlyeq}

\newcommand{\twoheadrightarrowtail}{\twoheadrightarrow \hspace{-12.19pt} \rightarrowtail}

\newcommand{\mysub}[1]{\hbox{\smaller[3]{$#1$}}}

\newcommand{\sub}[2]{#1_{\hbox{\smaller[3]{$#2$}}}}

\newcommand{\subAlt}[2]{#1_{\raisebox{0.95 ex}[0 ex][0 ex]{\smaller[3]{$#2$}}}}

\newcommand{\subM}[1]{#1_{\mysub{M}}}

\newcommand{\outSigma}{\widehat{\Sigma}}

\newcommand{\outMathfrakI}{\widehat{\mathfrak{I}}}

\newcommand{\outMoutMathfrakI}{\widehat{\makebox*{$m$}{$M$}}_{\mysub{\outMathfrakI}}}

\let\leq=\leqslant
\let\le=\leqslant
\let\geq=\geqslant
\let\ge=\geqslant

\newcommand{\sT}{\mid}
\definecolor{dark-gray}{gray}{0.35}

\newcommand{\rkk}{\hbox{\sf rk}\xspace}

\def \no#1#2#3 {{\bf #1} (#3), #2.}
\def \eds#1#2 {#1, #2.}

\definecolor{grey}{rgb}{0.9, 0.9, 0.9}
\newlength{\savefboxrule}
\newtheorem{mydef}{Definition} %
\newtheorem{myconj}{Conjecture} %
\newtheorem{mylemma}{Lemma} %
\newtheorem{mytheorem}{Theorem} %
\newtheorem{mycorollary}{Corollary} %
\newtheorem{myremarks}{Remarks} %
\newtheorem{myremark}[myremarks]{Remark} %
\newtheorem{myproblem}{Problem} %

\newcommand{\Vars}[1]{\mathrm{Vars}(#1)}

\newcommand{\COMMENT}[1]{}

\newcommand{\disj}[2]{#1\cap #2=\emptyset}

\newcommand{\card}[1]{{\left|{#1}\right|}}

\newcommand{\defAs}{\coloneqq}

\newcommand{\Nats}{\mathbb{N}}

\newcommand{\true}{{\bf true}}
\newcommand{\Pow}{{\mathrm{pow}}}
\newcommand{\pow}[1]{\Pow({#1})}

\newcommand{\powAst}{\Pow^{\ast}}
\newcommand{\powAstot}{\Pow_{1,2}^{\ast}}
\newcommand{\powOt}{\Pow_{1,2}}

\newcommand{\powast}[1]{\powAst({#1})}
\newcommand{\powot}[1]{\powOt({#1})}
\newcommand{\powastot}[1]{\powAstot({#1})}

\newcommand{\Nodes}{\mathcal{N}}

\newcommand{\st}{\,\texttt{|}\:}

\newcommand{\rk}{\hbox{\sf rk}\;}
\newcommand{\dom}{\hbox{\sf dom}}
\newcommand{\Places}{\mathcal{P}}
\newcommand{\TARGETS}{\mathcal{T}}
\newcommand{\Targets}[1]{\TARGETS({#1})}

\newcommand{\myphi}{\Phi}
\newcommand{\mypsi}{\Psi}
\newcommand{\boldcalV}{\mbox{\boldmath$\mathcal{V}$}}

\newcommand{\MLSSP}{\textnormal{\textsf{MLSSP}}\xspace}
\newcommand{\MLSSPF}{\textnormal{\textsf{MLSSPF}}\xspace}
\newcommand{\HF}{\textnormal{\textsf{HF}}\xspace}
\newcommand{\MLS}{\textnormal{\textsf{MLS}}\xspace}

\newcommand{\MLSC}{\hbox{\textnormal{\textsf{MLS}}\hspace{-1.3pt}\raisebox{.9pt}{$\times$}}\xspace}
\newcommand{\MLSuC}{\hbox{\textnormal{\textsf{MLS}}\hspace{-.5pt}\raisebox{.9pt}{$\otimes$}}\xspace}

\newcommand{\BST}{\textnormal{\textsf{BST}}}
\newcommand{\BSTuCsub}{\textnormal{\textsf{BST}}\xspace\raisebox{.9pt}{$\otimes_{_{\subseteq}}$}\xspace}
\newcommand{\BSTC}{\hbox{\textnormal{\textsf{BST}}\hspace{-1.3pt}\raisebox{.9pt}{$\times$}}\xspace}
\newcommand{\BSTuC}{\hbox{\textnormal{\textsf{BST}}\hspace{-.5pt}\raisebox{.9pt}{$\otimes$}}\xspace}

\newcommand{\BS}{\textnormal{\textsf{BST}}\xspace\raisebox{.9pt}\xspace}

\newcommand{\NP}{\textnormal{\textsf{NP}}\xspace}
\newcommand{\HTP}{\textnormal{\textsf{HTP}}\xspace}

\newcounter{instr}

\newcounter{instrb}
\newcommand{\ninstrb}{\refstepcounter{instrb}\textcolor{dark-gray}{\footnotesize{\theinstrb.}} \'}
\newcommand{\commentout}[1]{}

\begin{document}

\title{Decidability of the satisfiability problem for\\ Boolean set theory with the unordered\\ Cartesian product operator\thanks{We gratefully acknowledge partial support from the projects STORAGE and MEGABIT -- Universit\`{a} degli Studi di Catania, PIAno di inCEntivi per la RIcerca di Ateneo 2020/2022 (PIACERI), Linea di intervento 2.}}
\author{Domenico Cantone and  Pietro Ursino \\
\emph{Dipartimento di Matematica e Informatica, Universit\`a di Catania}\\
\emph{Viale Andrea Doria 6, I-95125 Catania, Italy.}  \\
\mbox{E-mail:} \texttt{domenico.cantone@unict.it,pietro.ursino@unict.it}
}

\maketitle
\begin{abstract}
The satisfiability problem for multilevel syllogistic extended with the Cartesian product operator (\MLSC) is a long-standing open problem in computable set theory. For long, it was not excluded that such a problem were undecidable, due to its remarkable resemblance with the well-celebrated Hilbert's tenth problem, as it was deemed reasonable that union of disjoint sets and Cartesian product might somehow play the roles of integer addition and multiplication.

To dispense with nonessential technical difficulties, we report here about a positive solution to the satisfiability problem for a slight simplified variant of \MLSC, yet fully representative of the combinatorial complications due to the presence of the Cartesian product, in which membership is not present and the Cartesian product operator is replaced with its unordered variant.

We are very confident that such decidability result can be generalized to full \MLSC, though at the cost of considerable technicalities. 

\end{abstract}

\section*{Introduction}
The decision problem in set theory has been studied quite thoroughly in the last decades, giving rise to the field of Computable Set Theory \cite{CFO89}. The initial goal was the mechanical formalisation of mathematics with a proof verifier based on the set-theoretic formalism \cite{OS02,COSU03,OCPS06,SchCanOmo11}, but soon a foundational interest aimed at the identification of the boundary in set theory between the decidable and the undecidable became more and more compelling.

The precursor fragment of set theory investigated for decidability was \MLS, which stands for Multi-Level Syllogistic. \MLS consists of the quantifier-free formulae of set theory involving only the Boolean set operators $\cup$, $\cap$, $\setminus$ and the relators $=$ and $\in$, besides set variables (assumed to be existentially quantified). The satisfiability problem (s.p., briefly) for  \MLS has been solved in the seminal paper \cite{FOS80a}, and its \NP-completeness has later been proved in \cite{COP90}. Following that, several extensions of \MLS with various combinations of the set operators $\{\cdot\}$ (singleton), $\Pow$ (power set), $\medcup$ (unary union), $\medcap$ (unary intersection), $\rkk$ (rank), etc., and of the set predicates rank comparison, cardinality comparison, finiteness, etc., have been also proved decidable over the years.\footnote{The monographs \cite{CFO89,COP01,SchCanOmo11,OPT17,CU18} provide a rather comprehensive account.}

However, the s.p.\ for the extension \MLSC of \MLS with the Cartesian product $\times$,\footnote{For definiteness, we may assume that the Cartesian product is expressed in terms of Kuratowski's ordered pairs $(u,v) \defAs \{\{u\}, \{u,v\}\}$.} proposed by the first author since the middle 80s, soon appeared to be very challenging and resisted several efforts to find a solution, either positive or negative. As a matter of fact, for long it was not excluded that the s.p.\ for \MLSC were undecidable (in particular, when restricted to finite models), due to its remarkable resemblance with the well-celebrated Hilbert's Tenth problem (\HTP, for short), posed by David Hilbert at the beginning of last century \cite{Hilbert-02}.\footnote{We recall that \HTP asks for a uniform procedure that can determine in a finite number of steps whether any given Diophantine polynomial equation with integral coefficients is solvable in integers. In 1970, it was shown that no algorithmic procedure exists for \HTP, as result of the combined efforts of M.\ Davis, H.\ Putnam, J.\ Robinson, and Y.\ Matiyasevich (DPRM theorem, see \cite{Rob,DPR61,Mat70}).} Indeed, it was deemed reasonable that the union of disjoint sets and the Cartesian product might somehow play the roles of integer addition and multiplication in \HTP, respectively, in consideration of the fact that $|s \cup t| = |s| + |t|$, for any disjoint sets $s$ and $t$, and $|s \times t| = |s| \cdot |t|$, for any sets $s$ and $t$.

Attempts to solve the s.p.\ for \MLSC helped shaping the development of computable set theory and led to the introduction of the powerful technique of \emph{formative processes},\footnote{See \cite{CU18} for a quite friendly introduction.} which has been at the base of the highly technical solutions to the decision problems for the extension \MLSSP of \MLS with the power set and the singleton operators \cite{COU02} and the extension \MLSSPF with the finiteness predicate too \cite{CU14}.  

\smallskip

In this paper, we provide a positive solution to the s.p., both unrestricted and restricted to (hereditarily) finite models, for the fragment of set theory dubbed \BSTuC, which is closely related to \MLSC. The fragment \BSTuC (which stands for Boolean Set Theory with the \emph{unordered Cartesian product} $\otimes$) is obtained by dropping the membership predicate $\in$ from \MLSC and by replacing the (ordered) Cartesian product operator $\times$ with its unordered variant $\otimes$, where $s \otimes t$ is the collection of all unordered pairs $\{u,v\}$ such that $u \in s$ and $v \in t$, namely $s \otimes t \coloneqq \big\{ \{u,v\} \st u \in s \wedge v \in t \big\}$ (for any sets $s$ and $t$).

Notice that none of the above two changes affects the aforementioned resemblance with \HTP. The reason why we chose to address here the case of \BSTuC rather than the one of \MLSC is that in doing so we can get rid of irrelevant features that would only make our analysis much more technical. Nevertheless, we plan to report on the s.p.\ for the fragments \BSTC and \MLSC in a future paper.\footnote{Naturally, \BSTC is Boolean Set Theory with Cartesian product.}

For both variants of the s.p., we shall provide  nondeterministic exponential decision procedures. These will be expressed in terms of the existence of a special graph, called \emph{$\otimes$-graph}, enjoying a certain connectivity property of \emph{accessibility}. Given a \BSTuC-formula $\myphi$ to be tested for satisfiability, in the case of the ordinary s.p.\ it will be enough to require that a candidate accessible $\otimes$-graph \emph{fulfills} $\myphi$, whereas in the case of the finite s.p.\ it will be additionally required that the $\otimes$-graph admits also a kind of topological order. In both cases, it will be shown that, when satisfied, these conditions (which are also necessary) ensure that the $\otimes$-graph can be used as a kind of flow graph to build a model for $\myphi$  in denumerably many steps (in the case of the ordinary s.p.) or in a bounded finite number of steps (in the case of the (hereditarily) finite s.p.). Such construction process is a simplified form of the formative processes mentioned before.

\smallskip

As shown in \cite{Schw78,CCS90}, the finite s.p.\ for the extension of \MLS with cardinality comparison, namely the the two-place predicate $|\cdot| \leq |\cdot|$ for cardinality comparison, where $|s| \leq |t|$ holds if and only if the cardinality of $s$ does not exceed that of $t$, can be reduced to purely existential Presburger arithmetic, which is known to be \NP-complete (see \cite{Sca}).
On the other hand, when \BSTC or \BSTuC is enriched with cardinality comparison,  the s.p.\ for the resulting extension become undecidable, since \HTP would be reducible to it, much as proved in \cite{CCP90} and \cite{COP20} for \MLSuC. This is clear evidence that the decision problem for both \BSTC and \BSTuC is very close to the border of decidability.

\smallskip
\centerline{---------------------------------------------}
\smallskip

The paper is organized as follows. In Section~\ref{se:BSTuC} we introduce the fragment of our interest \BSTuC through its syntax and semantics. In particular, semantics is presented in terms of satisfying partitions, and it is shown that such approach leads easily to the decidability of the purely Boolean subset of \BSTuC. It is also defined a useful variant of the \emph{intersecting power set} operator, in terms of which the unordered Cartesian product is easily expressible. Subsequently, in Section~\ref{se:ordinary satisfiability problem}, we introduce the central notion of accessible $\otimes$-graphs, together with that of fullfilment of a \BSTuC-formula by an accessible $\otimes$-graph, and we prove that any satisfiable \BSTuC-formula is fulfilled by a suitable accessible $\otimes$-graph. We also prove that such condition is sufficient for the satisfiability of $\myphi$, by describing in details a construction process that uses an accessible $\otimes$-graph fulfilling $\myphi$  as a kind of flow graph to build a model for $\myphi$ in denumerably many steps. Afterward, in Section~\ref{se:finite satisfiability}, we introduce the notion of ordered $\otimes$-graphs, and prove that the existence of an ordered $\otimes$-graph fulfilling a given \BSTuC-formula $\myphi$ is a necessary and sufficient condition for $\myphi$ to be (hereditarily) finitely satisfiable. Finally, in Section~\ref{se:concluding remarks}, we discuss some plans for future research.

\vspace{1cm}

\section{The fragment \BSTuC}\label{se:BSTuC}
$\BSTuC$ is the quantifier-free fragment of set theory consisting of the propositional closure of atoms of the following types:
\[
x=y \cup z \/,  \quad x=y \cap z\/,  \quad x=y \setminus z\/,  \quad x = y \otimes z \/,   \quad  x\subseteq y\/
\]
where $x,y,z$ stand for set variables. For any \BSTuC-formula $\myphi$, we denote by $\Vars{\myphi}$ the collection of set variables occurring in it.

\subsection{Semantics of \BSTuC}
The semantics of \BSTuC is defined in a very natural way in terms of \emph{set assignments}.  

A \textsc{set assignment} $M$ is any map from a collection $V$ of set variables (called the \textsc{variables domain of $M$} and denoted $\dom(M)$) into the von Neumann universe $\boldcalV$ of all well-founded sets.

We recall that $\boldcalV$ is a cumulative hierarchy constructed in stages by transfinite recursion over the class $\mathit{On}$ of all ordinals. Specifically, $\boldcalV \defAs \bigcup_{\alpha \in \mathit{On}} \mathcal{V}_{\alpha}$ where, recursively, $\mathcal{V}_{\alpha} \defAs \bigcup_{\beta<\alpha} \pow{\mathcal{V}_{\beta}}$, for every $\alpha \in \mathit{On}$, with $\pow{\cdot}$ denoting the powerset operator. Based on such construction, we can readily define the \textsc{rank} of any well-founded set $s \in \boldcalV$, denoted $\rk{s}$, as the least ordinal $\alpha$ such that $s \subseteq \mathcal{V}_{\alpha}$. The collection of the sets of finite rank, hence belonging to $\mathcal{V}_{\alpha}$ for some finite ordinal $\alpha$, forms the set \HF of the \textsc{hereditarily finite sets}. Thus, $\HF = \mathcal{V}_{\omega}$, where $\omega$ is the first limit ordinal, namely the smallest non-null ordinal with no immediate predecessor.

Given a set assignment $M$ and a collection of variables $W \subseteq \dom(M)$, we put $MW \defAs \{Mv \sT v \in W\}$. The \textsc{set domain of $M$} is defined as the set $\bigcup MV = \bigcup_{v \in V}Mv$. The \textsc{rank of $M$} is the rank of its set domain, namely, $\rk M  \defAs  \rkk (\bigcup MV)$ (so that, when $V$ is finite, $\rk M = \max_{v \in V} ~\rk Mv$).
A set assignment $M$ is \textsc{finite} (resp., \textsc{hereditarily finite}), if so is its set domain.

\smallskip

Operators and relators of \BSTuC are interpreted according to their usual semantics. Thus, given a set assignment $M$, for any $x,y,z \in \dom(M)$ we put:
\begin{align*}
M(x \star y) \defAs  Mx \star My, 
\end{align*}
where $\star \in \{\cup,\cap,\setminus\}$, and 
\begin{align*}
M(x = y \star z) = \text{\bf true} &\qquad \xleftrightarrow{\hbox{\smaller[5]{$~\mathit{Def}$~}}} \qquad Mx = M(y \star z),\\
M(x \subseteq y) = \text{\bf true} &\qquad\xleftrightarrow{\hbox{\smaller[5]{$~\mathit{Def}$~}}} \qquad Mx \subseteq My.
\end{align*}
Finally, we put recursively
\[
\begin{aligned}
M(\neg \myphi) &\defAs \neg M \myphi, \quad & M(\myphi \wedge \mypsi) &\defAs  M \myphi \wedge M \mypsi,\\
M(\myphi \vee \mypsi) &\defAs  M \myphi \vee M \mypsi, \quad & M(\myphi \rightarrow \mypsi) &\defAs  M \myphi \rightarrow M \mypsi, & \quad \text{etc.,}
\end{aligned}
\]
for all \BSTuC-formulae $\myphi$ and $\mypsi$ such that $\Vars{\myphi}, \Vars{\mypsi} \subseteq \dom(M)$.

Given a \BSTuC-formula $\myphi$, a set assignment $M$ over $\Vars{\myphi}$ is said to \textsc{satisfy} $\myphi$ if $M \Phi = \true$ holds, in which case we also write $M \models \myphi$ and say that $M$ is a \textsc{model} for $\myphi$. If $\myphi$ has a model, we say that $\myphi$ is \textsc{satisfiable}; otherwise, we say that $\myphi$ is \textsc{unsatisfiable}. If $M \models \myphi$ and $M$ is finite (resp., hereditarily finite), then $\myphi$ is \textsc{finitely satisfiable} (resp., \textsc{hereditarily finitely satisfiable}).

Two \BSTuC-formulae $\myphi$ and $\mypsi$ are \textsc{equisatisfiable} when $\myphi$ is satisfiable if and only if so is $\mypsi$, possibly by distinct models.

The \textsc{decision problem} or \textsc{satisfiability problem} for \BSTuC is the problem of establishing algorithmically whether any given \BSTuC-formula is satisfiable or not by some set assignment. 

By restricting to (hereditarily) finite set assignments, one can define in the obvious way the \textsc{(hereditarily) finite satisfiability problem} for \BSTuC.

\subsection{Satisfiability by partitions}

A \textsc{partition} is a collection of pairwise disjoint non-null sets, called the \textsc{blocks} of the partition. The union $\bigcup \Sigma$ of a partition $\Sigma$ is its \textsc{domain}.

Let $V$ be a finite collection of set variables and $\Sigma$ a  partition. Also, let $\mathfrak{I} \colon V \rightarrow \pow{\Sigma}$ be any map. In a very natural way, the map $\mathfrak{I}$ induces a set assignment $M_{\mysub{\mathfrak{I}}}$ over $V$ definded by:
\[
\textstyle
M_{\mysub{\mathfrak{I}}} v \defAs \bigcup \mathfrak{I}(v)\/, \qquad \text{for $v \in V$\/.}
\]
We refer to the map $\mathfrak{I}$ (or to the pair $(\Sigma, \mathfrak{I})$, when we want to emphasize the partition $\Sigma$) as a \textsc{partition assignment}.

\begin{mydef}\label{def:satisfiability}\rm
Let $\Sigma$ be a partition and $\mathfrak{I} \colon V \rightarrow \pow{\Sigma}$ be a partition assignment over a finite collection $V$ of set variables. Given a \BSTuC-formula $\myphi$ such that $\Vars{\myphi} \subseteq V$, we say that $\mathfrak{I}$ \textsc{satisfies $\myphi$}, and write $\mathfrak{I} \models \myphi$, when the set assignment $M_{\mysub{\mathfrak{I}}}$ induced by $\mathfrak{I}$ satisfies $\myphi$ (equivalently, one may say that $\Sigma$ \textsc{satisfies $\myphi$ via the map $\mathfrak{I}$}, and write $\Sigma/\mathfrak{I} \models \myphi$, if we want to emphasize the partition $\Sigma$). We say that $\Sigma$ \textsc{satisfies} $\myphi$, and write $\Sigma \models \myphi$, if $\Sigma$ satisfies $\myphi$ via some map $\mathfrak{I} \colon  V \rightarrow \pow{\Sigma}$.
\end{mydef}

The following result can be proved immediately.
\begin{mylemma}\label{wasB}
If a \BSTuC-formula is satisfied by a partition $\Sigma$, then it is satisfied by any partition $\overline \Sigma$ that includes $\Sigma$ as a subset, namely such that $\Sigma \subseteq \overline \Sigma$.
\end{mylemma}

Plainly, a \BSTuC-formula $\myphi$ satisfied by some partition is satisfied by a set assignment. Indeed, if $\Sigma \models \myphi$, then $\Sigma/\mathfrak{I} \models \myphi$ for some map $\mathfrak{I} \colon  V \rightarrow \pow{\Sigma}$, and therefore $M_{\mysub{\mathfrak{I}}} \models \myphi$. The converse holds too. In fact, let us assume that $M \models \myphi$, for some set\index{set} assignment $M$ over the collection $V = \Vars{\myphi}$ of the set variables occurring in $\myphi$, and let $\subM{\Sigma}$ be the \textsc{Venn partition} induced by $M$, namely
\[
\sub{\Sigma}{M} \defAs \Big\{ \medcap MV' \setminus \medcup M(V \setminus V') \st \emptyset \neq V' \subseteq V \Big\} \setminus \big\{ \,\emptyset\,\big\}.\footnotemark
\]
\footnotetext{Hence, we have:
\begin{enumerate}[label=-]
\item $(\forall \sigma \in \sub{\Sigma}{M})(\forall v \in V)(\sigma \cap Mv = \emptyset \vee \sigma \subseteq Mv)$,

\item $(\forall \sigma,\sigma' \in \sub{\Sigma}{M}) \big( (\forall v \in V) (\sigma \subseteq Mv \leftrightarrow \sigma' \subseteq Mv) \leftrightarrow \sigma = \sigma' \big)$, and

\item $\medcup \Sigma = \medcup MV$.
\end{enumerate}
}
Let $\subM{\mathfrak{I}} \colon V \rightarrow \pow{\subM{\Sigma}}$ be the map defined by 
\[
\subM{\mathfrak{I}}(v) \defAs \{ \sigma \in \subM{\Sigma} \st \sigma \subseteq Mv\}\/, \qquad \text{for $v \in V$.}
\]
It is an easy matter to check that the set assignment induced by $\subM{\mathfrak{I}}$ is just $M$. Thus $\subM{\Sigma}/\subM{\mathfrak{I}} \models \myphi$, and therefore $\subM{\Sigma} \models \myphi$, proving that $\myphi$ is satisfied by some partition, in fact by the Venn partition induced by $M$, whose size is at most $2^{|V|}-1$.

Thus, the notion of satisfiability by set assignments and that of satisfiability by partitions coincide.

As a by-product of Lemma~\ref{wasB} and the above considerations, we also have:
\begin{mylemma}\label{wasA}
Every \BSTuC-formula $\myphi$ with $n$ distinct variables is satisfiable if and only if it is satisfied by some partition with $2^{n}-1$ blocks.
\end{mylemma}

\subsection{Normalization of \BSTuC-formulae}

By applying disjoint normal form and the simplification rules illustrated in \cite{CU18}, the satisfiability problem for \BSTuC can be reduced to the satisfiability problem for \textsc{normalized conjunctions} of \BSTuC, namely
conjunctions of \BSTuC-literals of the following restricted types:%
\begin{gather}\label{formula}
  x=y \cup z \/,\ \  x=y \setminus z\/,\ \    x = y \otimes z\/,\ \  x \neq y\/,
  \end{gather}
where $x,y,z$ stand for set variables. Indeed, it is enough to observe that:
\begin{itemize}
\item[--] $x \subseteq y$ is equivalent to $x = x \cap y$;

\item[--] the terms $y \cap z$ and $y \setminus (y \setminus z)$ are equivalent, so  an atom of the form $x = y \cap z$ is equisatisfiable with the conjunction $x = y \setminus y' \wedge y' = y \setminus z$, where $y'$ stands for any fresh set variable;

\item[--] each negative literal of the form $x \neq y \star z$ (with $\star \in \{\cup, \cap,\setminus,\otimes\}$) is equisatisfiable with the conjunction $x' = y \star z \wedge x' \neq x$, where $x'$ stands for any fresh set variable.
\end{itemize}

\subsection{The Boolean case}

In the restricted case of \textsc{$\BST$-conjunctions}, namely conjunctions of Boolean literals of the form 
\begin{equation*}\label{Boolean literals}
x=y \cup z, \quad x=y \setminus z, \quad x \neq y,
\end{equation*}
the satisfiability status by a given partition $\Sigma$ does not depend in any way on the internal structure of its blocks, but just on their numerousness. This is proved in the following lemmas.

\begin{mylemma}\label{partitionAssignmentBoolean}
Let $\Sigma$ be a partition and let $\mathfrak{I} \colon V \rightarrow \pow{\Sigma}$ be a partition assignment over a (finite) set of variables $V$. Then, for all $x,y,z \in V$ and $\star \in \{\cup,\setminus\}$, we have:
\begin{enumerate}[label=(\alph*)]
\item\label{partitionAssignmentBooleanA} $\mathfrak{I} \models x=y \star z \quad \Longleftrightarrow \quad \mathfrak{I}(x) = \mathfrak{I}(y) \star \mathfrak{I}(z)$,

\item\label{partitionAssignmentBooleanB} $\mathfrak{I} \models x\neq y\quad \Longleftrightarrow \quad \mathfrak{I}(x) \neq \mathfrak{I}(y)$.
\end{enumerate}
\end{mylemma}
\begin{proof}\belowdisplayskip=-13pt%
It is enough to observe that since $\Sigma$ is a partition (and therefore its blocks are nonempty and mutually disjoint), then for all $x,y,z \in V$ and $\star \in \{\cup, \setminus\,\}$ we have:
\begin{align*}
\mathfrak{I} \models x=y \star z \quad &\Longleftrightarrow \quad \medcup \mathfrak{I}(x) = \medcup \mathfrak{I}(y) \star \medcup\mathfrak{I}(z) \\
&\Longleftrightarrow \quad \medcup \mathfrak{I}(x) = \medcup \big(\mathfrak{I}(y) \star \mathfrak{I}(z)\big) \\
&\Longleftrightarrow \quad \mathfrak{I}(x) = \mathfrak{I}(y) \star \mathfrak{I}(z) \\ %
\shortintertext{and}
\mathfrak{I} \models x\neq y \quad &\Longleftrightarrow \quad \medcup \mathfrak{I}(x) \neq \medcup \mathfrak{I}(y) \\
&\Longleftrightarrow \quad \mathfrak{I}(x) \neq \mathfrak{I}(y). %
\end{align*}
\end{proof}

Satisfiability of $\BST$-conjunction can be expressed in purely combinatorial terms by means of \emph{fulfilling maps}.
\begin{mydef}\label{Fulfilling-BST}
\rm Let $\myphi$ be a $\BST$-conjunction, and let $\mathfrak{F}\colon V \rightarrow \pow{{\Pow^{+}(V)}}$ be any map, where $V \defAs \Vars{\myphi}$ and $\Pow^{+}(V) \defAs \pow{V}\setminus \emptyset$. We say that the map $\mathfrak{F}$ \textsc{fulfills} $\myphi$ provided that:

\begin{enumerate}[label=(\alph*)]
\item\label{Fulfilling-BSTa} $\mathfrak{F}(x) = \mathfrak{F}(y) \star \mathfrak{F}(z)$, for each conjunct $x = y \star z$ in $\myphi$, with $\star \in \{\cup,\setminus\}$;

\item\label{Fulfilling-BSTb} $\mathfrak{F}(x) \neq \mathfrak{F}(y)$, for each conjunct $x \neq y$ in $\myphi$.
\end{enumerate}
A map $\mathfrak{F}$ satisfying conditions \ref{Fulfilling-BSTa} and \ref{Fulfilling-BSTb} above will be called a \textsc{fulfilling map} for $\myphi$.
\end{mydef}

\begin{myremark}\rm
In Section~\ref{Fulfillment by accessible graphs}, fulfilling maps will be defined differently, in the context of the s.p.\ for \BSTuC. However, such an overloading should create no problems.
\end{myremark}

\begin{mylemma}\label{le:BST-fulfillment}
A $\BST$-conjunction is satisfiable if and only if it is fulfilled by some fulfilling map.
\end{mylemma}
\begin{proof}
Let $\myphi$ be a $\BST$-conjunction and let $V \defAs \Vars{\myphi}$.

For the sufficiency part, let us first assume that $\myphi$ is satisfied by a partition $\Sigma$ via a certain map $\mathfrak{I}\colon V \rightarrow \pow{\Sigma}$. For each $\sigma \in \Sigma$, we put:
\[
V_{\sigma} \defAs \{v \in V \st \sigma \in \mathfrak{I}(v)\}.
\]
Then, we define the map $\sub{\mathfrak{F}}{\mathfrak{I}}\colon V \rightarrow \pow{{\Pow^{+}(V)}}$ by putting
\[
\sub{\mathfrak{F}}{\mathfrak{I}}(x) \defAs \{V_{\sigma} \st \sigma \in \mathfrak{I}(x)\}, \qquad \text{for } x \in V.
\]
Preliminarily, we observe that 
\begin{equation}\label{prelimEquation}
V_{\sigma} \in \sub{\mathfrak{F}}{\mathfrak{I}}(x) ~\longleftrightarrow~ \sigma \in \mathfrak{I}(x),\qquad \text{for $\sigma \in \Sigma$ and $x \in V$}.
\end{equation}
Indeed, if $V_{\sigma} \in \sub{\mathfrak{F}}{\mathfrak{I}}(x)$, then $V_{\sigma} = V_{\sigma'}$, for some $\sigma' \in \mathfrak{I}(x)$. But then, since $\sigma \in \mathfrak{I}(v) \longleftrightarrow \sigma' \in \mathfrak{I}(v)$, for all $v \in V$, we have $\sigma \in \mathfrak{I}(x)$.

Thus, for every conjunct $x = y \star z$ in $\myphi$ (with $\star \in \{\cup,\setminus\}$), we have:
\begin{align*}
\sub{\mathfrak{F}}{\mathfrak{I}}(x) &= \{V_{\sigma} \st \sigma \in \mathfrak{I}(x)\}\\
&= \{V_{\sigma} \st \sigma \in \mathfrak{I}(y) \star \mathfrak{I}(z)\} &&\text{(by Lemma~\ref{partitionAssignmentBoolean}\ref{partitionAssignmentBooleanA})}\\
&= \{V_{\sigma} \st \sigma \in \mathfrak{I}(y)\} \star \{V_{\sigma} \st \sigma \in \mathfrak{I}(z)\} &&\text{(by \eqref{prelimEquation})}\\
&= \sub{\mathfrak{F}}{\mathfrak{I}}(y) \star \sub{\mathfrak{F}}{\mathfrak{I}}(z).
\end{align*}
Likewise, for every conjunt $x \neq y$ in $\myphi$, by Lemma~\ref{partitionAssignmentBoolean}\ref{partitionAssignmentBooleanB} and \eqref{prelimEquation}, we have:
\[
\sub{\mathfrak{F}}{\mathfrak{I}}(x) = \{V_{\sigma} \st \sigma \in \mathfrak{I}(x)\} \neq \{V_{\sigma} \st \sigma \in \mathfrak{I}(y)\} = \sub{\mathfrak{F}}{\mathfrak{I}}(y).
\]
Thus, the map $\sub{\mathfrak{F}}{\mathfrak{I}}$ fulfills $\myphi$.

\smallskip

Conversely, for the necessity part, let us assume that $\myphi$ is fulfilled by a map $\mathfrak{F}\colon V \rightarrow \pow{{\Pow^{+}(V)}}$. Let $\sub{\Sigma}{\mathfrak{F}}$ be any partition of size $|\Pow^{+}(V)|$, and let $\beta \colon \Pow^{+}(V) \twoheadrightarrowtail \sub{\Sigma}{\mathfrak{F}}$ be any bijection from $\Pow^{+}(V)$ onto $\sub{\Sigma}{\mathfrak{F}}$. Let us define the map $\sub{\mathfrak{I}}{\mathfrak{F}} \colon V \rightarrow \pow{\sub{\Sigma}{\mathfrak{F}}}$ by setting
\[
\sub{\mathfrak{I}}{\mathfrak{F}}(x) \defAs \beta[\mathfrak{F}(x)], \qquad \text{for } x \in V.
\]
But then, for every literal $x = y \star z$ in $\myphi$ (with $\star \in \{\cup,\setminus\}$), we have:
\begin{align*}
\sub{\mathfrak{I}}{\mathfrak{F}}(x) &= \beta[\mathfrak{F}(x)] \\
&= \beta[\mathfrak{F}(y) \star \mathfrak{F}(z)]\\
&= \beta[\mathfrak{F}(y)] \star \beta[\mathfrak{F}(z)]\\
&= \sub{\mathfrak{I}}{\mathfrak{F}}(y) \star \sub{\mathfrak{I}}{\mathfrak{F}}(z).
\end{align*}
Hence, by Lemma~\ref{partitionAssignmentBoolean}\ref{partitionAssignmentBooleanA}, $\sub{\mathfrak{I}}{\mathfrak{F}} \models x = y \star z$.

Similarly, for every conjunct $x \neq y$ in $\myphi$, we have
\[
\sub{\mathfrak{I}}{\mathfrak{F}}(x) = \beta[\mathfrak{F}(x)] \neq \beta[\mathfrak{F}(y)] = \sub{\mathfrak{I}}{\mathfrak{F}}(y),
\]
proving that, by Lemma~\ref{partitionAssignmentBoolean}\ref{partitionAssignmentBooleanB}, $\sub{\mathfrak{I}}{\mathfrak{F}} \models x \neq y$.

Hence, in conclusion, we have $\sub{\mathfrak{I}}{\mathfrak{F}} \models \myphi$, proving that $\myphi$ is satisfiable.
\end{proof}

An immediate consequence of the preceding lemma is the following result.

\begin{mylemma}\label{le:every}
A $\BST$-conjunction with $n$ distinct variables is satisfiable if and only if it is satisfied by \emph{every} partition of size $2^{n}-1$.
\end{mylemma}

In fact, in the light of \cite[Lemma 2.36, p.\ 42]{CU18} and Lemma~\ref{wasA}, one can prove the following stronger result:
\begin{mylemma}\label{was_Lemma5}
A $\BST$-conjunction involving $n$ distinct variables is satisfiable if and only if it is satisfied by \emph{every} partition of size $n - 1$.
\end{mylemma}

In view of the preceding lemma, Lemma~\ref{le:BST-fulfillment} can be so strengthened:
\begin{mylemma}\label{le:small fulfilling map}
A \BST-conjunction over a set $V$ of variables is satisfiable if and only if it is fulfilled by a map $\mathfrak{F}\colon V \rightarrow \pow{{\Pow^{+}(V)}}$ such that $|\medcup\mathfrak{F}[V]| \leq |V|-1$.
\end{mylemma}

Both previous two lemmas readily yield that the satisfiability problem for $\BST$-conjunctions can be solved in nondeterministic polynomial time, namely it belongs to the class \NP.

As shown in \cite{CDMO19}, the satisfiability problem for conjunctions of Boolean literals of the form $t_{1} \neq t_{2}$, where $t_{1}$ and $t_{2}$ are set terms involving only variables and the set difference operator `$\setminus$', is \NP-complete. Therefore, we have:
\begin{mylemma}
The satisfiability problem for $\BST$-conjunctions is \emph{\NP}-complete.
\end{mylemma}

\subsection{Dealing also with literals of type $x = y \otimes z$: the intersecting power set operator}

We shall express the conditions that take also care of literals in \BSTuC of  type $x = y \otimes z$ by means of some useful variants of the power set operator and the \textsc{intersecting power set operator} $\powAst$, defined by
\begin{align*}
\powast{S} &\defAs \big\{t \subseteq \medcup S \st t \cap s \neq \emptyset, \text{ for every } s \in S \big\}
\end{align*}
and introduced in \cite{Can91} in connection with the solution of the satisfiability problem for a fragment of set theory involving the power set  and the singleton operators.\footnote{Several properties of the operator $\powAst$ are listed in \cite[pp.\ 16--20]{CU18}.} Specifically, for any set $S$, we put:
\begin{align*}
\powastot{S} &\defAs \big\{t \in \powast{S} \st 1 \le |t| \le 2 \big\},\\
\powot{S} &\defAs \big\{t \in \pow{S} \st 1 \le |t| \le 2 \big\}.
\end{align*}
Thus,
\begin{enumerate}[label=-]
\item $\powastot{S}$ is the set of all the members of $\powast{S}$ of cardinality equal to either $1$ or $2$;
\item $\powot{S}$ is the collection of all the subsets of $S$ of cardinality equal to either $1$ or $2$.
\end{enumerate}

\medskip

We state next a useful injectivity property of the operator $\powAstot$.
\begin{mylemma}\label{injectivity-powAstot}
Let $\Sigma$ be a partition.
For all $\mathcal{B},\mathcal{B'} \subseteq \powot{\Sigma}$, we have
\[
\medcup \powAstot[\mathcal{B}] = \medcup \powAstot[\mathcal{B'}]  \quad \Longleftrightarrow \quad \mathcal{B} = \mathcal{B'}.
\]
\end{mylemma}

The $\powAstot$ operator is strictly connected with the unordered Cartesian operator $\otimes$, as we show next.

\begin{mylemma}\label{powastotAndOtimes}
For all sets $s$ and $t$ (not necessarily distinct), we have
\[
\powastot{\{s,t\}} = s \otimes t.
\]
\end{mylemma}
\begin{proof}
Plainly, $s \otimes t \subseteq \powastot{\{s,t\}}$. Indeed, if $u \in s \otimes t$, then
\[
1 \le |u| \le 2, \quad u \subseteq s \cup t, \quad \text{and} \quad u \cap s \neq \emptyset \neq u \cap t,
\]
so that $u \in \powastot{\{s,t\}}$.

Conversely, let $\{u,v\} \in \powastot{\{s,t\}}$. Then
\[
\{u,v\} \subseteq s \cup t \quad \text{and} \quad \{u,v\} \cap s \neq \emptyset \neq \{u,v\} \cap t.
\]
Without loss of generality, let us assume that $u \in s$. If $v \in t$, we are done. Otherwise, if $v \notin t$, then $v \in s$ (since $\{u,v\} \subseteq s \cup t$) and $u \in t$ (since $\{u,v\} \cap t \neq \emptyset$). Hence, $\{u,v\} \in s \otimes t$, proving that also the inverse inclusion $\powastot{\{s,t\}} \subseteq s \otimes t$ holds.
\end{proof}

The unordered Cartesian operator $\otimes$ enjoys the following distributive property.

\begin{mylemma}[Distributivity]\label{wasD}
For all sets $S$ and $T$, the following identity holds:
\begin{equation}\label{eq:wasD}
\medcup S \otimes \medcup T = \medcup \{ s \otimes t \emph{\st} s \in S,~ t \in T\}.
\end{equation}
\end{mylemma}

\begin{proof}
Let $u \in \medcup S \otimes \medcup T$. Then $u = \{u',u''\}$ for some $u' \in \medcup S$ and $u'' \in \medcup T$. Hence, $u' \in \overline s$ and $u'' \in \overline t$ for some $\overline s \in S$ and $\overline t \in T$, so that 
$
u = \{u',u''\} \in \medcup \{s \otimes t \st s\in S,~ t\in T\}.
$
Thus, 
\begin{equation}\label{eq:wasD-firstHalf}
\medcup S \otimes \medcup T \subseteq \medcup \{s \otimes t \st s\in S,~ t\in T\}.
\end{equation}

For the converse inclusion, let $u \in \medcup \{s \otimes t \st s\in S,~ t\in T\}$. Then $u \in \overline s \otimes \overline t$, for some $\overline s \in S$ and $\overline t \in T$, and therefore $u = \{u',u''\}$, for some $u' \in \overline s$ and $u'' \in \overline t$. Since $\overline s \subseteq \medcup S$ and $\overline t \subseteq \medcup T$, then $u \in \medcup S \otimes \medcup T$, and so
\[
\medcup \{s \otimes t \st s\in S,~ t\in T\} \subseteq \medcup S \otimes \medcup T.
\]
Together with the converse inclusion \eqref{eq:wasD-firstHalf}, the latter yields \eqref{eq:wasD}, completing the proof of the lemma.
\end{proof}

\section{The ordinary satisfiability problem for \BSTuC-conjunctions}\label{se:ordinary satisfiability problem}

A result like those contained in Lemmas~\ref{le:every} and \ref{was_Lemma5} cannot hold for \BSTuC, since the literals of type $x = y \otimes z$ force one to take into account also the internal structure of certain blocks in any partition that satisfies a given \BSTuC-conjunction. These are the \emph{$\otimes$-blocks}, which are defined next.

\begin{mydef}\label{otimes-blocks-upblocks}\rm
A subset $\Sigma^{*}$ of a partition $\Sigma$ is a \textsc{$\otimes$-subpartition} of $\Sigma$  if $\medcup \Sigma^{*} = \medcup \powAstot[\mathcal{B}]$, for some $\mathcal{B} \subseteq \Sigma \otimes \Sigma$.

We denote by $\sub{\Sigma}{\otimes}$ the $\subseteq$-maximal $\otimes$-subpartition of $\Sigma$ and we refer to its elements as the \textsc{$\otimes$-blocks} of $\Sigma$.\footnote{The definition of $\sub{\Sigma}{\otimes}$ is well given, since the collection of the $\otimes$-subpartitions of $\Sigma$ is closed under union.}
We also denote by $\sub{\Pi}{\otimes}$ the subset of $\Sigma \otimes \Sigma$ such that $\medcup \sub{\Sigma}{\otimes} = \medcup \powAstot[\sub{\Pi}{\otimes}]$ and we refer to its elements as \textsc{$\otimes$-upblocks} (`upblocks' for \emph{unordered pair of blocks}).\footnote{By Lemma~\ref{injectivity-powAstot}, the set $\sub{\Pi}{\otimes}$ is well defined.}
\end{mydef}

As a consequence of Lemma~\ref{le:every} (resp., Lemma~\ref{was_Lemma5}), to test whether a given conjunction $\myphi$ of $\BST$-conjunctions with $n$ distinct variables is satisfiable, it is enough to check whether a partition whatsoever with $2^{n}-1$ (resp., $n-1$) blocks satisfies $\myphi$.

In the case of \BSTuC-conjunctions with $n$ distinct variables, rather than checking a single partition for satisfiability, one would have to test a whole collection of doubly exponential size of partitions with $2^{n} - 1$ blocks. Remarkably, the partitions in such collection can be conveniently described by special graphs, called \emph{$\otimes$-graphs}, which enjoy a particular connectivity property termed \emph{accessibility}.

Given a \BSTuC-conjunction $\myphi$ to be tested for satisfiability, in the case of the ordinary s.p.\ it will be enough to find an accessible $\otimes$-graph which \emph{fulfills} $\myphi$, in the sense that will be soon made precise, whereas for the (hereditarily) finite s.p., besides accessibility and fulfillability, it will be additionally requested that the $\otimes$-graphs admit a ``weak'' topological order.

Both for the ordinary s.p.\ and for the (hereditarily) finite s.p.\, such an approach will yield nondeterministic exponential decision procedures in the number of distinct variables of the input formula.

\smallskip

Next, we provide precise definitions of the notions mentioned above. We begin with $\otimes$-graphs.

\subsection{$\otimes$-graphs}

\begin{mydef}[$\otimes$-graphs] \rm
A \textsc{$\otimes$-graph} $\mathcal{G}$ is a directed bipartite graph whose set of vertices comprises two disjoint parts: a set of \textsc{places} $\Places$ and a set of \textsc{nodes} $\Nodes$, where $\Nodes = \Places \otimes \Places$.\footnote{Thus, it is required that $\Places \cap (\Places \otimes \Places) = \emptyset$.} The edges issuing from each place $q$ are exactly all pairs $\langle q,B \rangle$ such that $q\in B \in \sub{\Nodes}{}$\/: these are called \textsc{membership edges}. The remaining edges of $\sub{\mathcal{G}}{}$, called \textsc{distribution} or \textsc{saturation edges}, go from nodes to places. When there is an edge $\langle B,q \rangle$ from a node $B$ to a place $q$, we say that $q$ is a \textsc{target} of $B$. The map $\sub{\TARGETS}{}$ over $\sub{\Nodes}{}$ defined by
\[
\sub{\TARGETS}{}(B) \defAs \{q \in \Places \st q \text{ is a target of } B\}, \quad \text{for } B \in \sub{\Nodes}{},
\]
is the \textsc{target map} of $\sub{\mathcal{G}}{}$. The \textsc{size} of $\sub{\mathcal{G}}{}$ is defined as the cardinality of its set of places $\sub{\Places}{}$. Plainly, a $\otimes$-graph $\sub{\mathcal{G}}{}$ is fully characterized by its target map $\sub{\TARGETS}{}$, since the sets of nodes and of places of $\sub{\mathcal{G}}{}$ are expressible as $\dom(\TARGETS)$ and $\medcup \dom(\TARGETS)$, respectively. When convenient, we shall explicitly write $\sub{\mathcal{G}}{} = (\Places, \Nodes, \TARGETS)$ for a $\otimes$-graph with set of places $\Places$, set of nodes $\Nodes$, and target map $\TARGETS$.
\end{mydef}

\smallskip

To better grasp the rationale behind the definition just stated of $\otimes$-graphs, it is helpful to illustrate how to construct the $\otimes$-graph $\sub{\mathcal{G}}{\Sigma}$ \textsc{induced by a given a partition} $\Sigma$. 

To begin with, we select a set of places $\sub{\Places}{\Sigma}$ of the same cardinality of $\Sigma$ such that $\sub{\Places}{\Sigma}$ and $\sub{\Places}{\Sigma} \otimes \sub{\Places}{\Sigma}$ are disjoint, and define the vertex set of $\sub{\mathcal{G}}{\Sigma}$ as the union $\sub{\Places}{\Sigma} \cup (\sub{\Places}{\Sigma} \otimes \sub{\Places}{\Sigma})$. The members of $\sub{\Places}{\Sigma} \otimes \sub{\Places}{\Sigma}$ (namely the nonempty subsets of $\sub{\Places}{\Sigma}$ having cardinality at most 2) will form the set of nodes $\sub{\Nodes}{\Sigma}$ of $\sub{\mathcal{G}}{\Sigma}$. Places in $\sub{\Places}{\Sigma}$ are intended to be an abstract representation of the blocks of $\Sigma$ via a bijection $q \mapsto q^{(\bullet)}$ from $\sub{\Places}{\Sigma}$ onto $\Sigma$. Likewise, nodes in $\sub{\Nodes}{\Sigma}$ are intended to represent the unordered Cartesian product of the blocks represented by their places. The disjoint sets $\sub{\Places}{\Sigma}$ and $\sub{\Nodes}{\Sigma}$ will form the parts of the bipartite graph $\sub{\mathcal{G}}{\Sigma}$ we are after. The bijection $(\bullet)$ can be naturally extended to nodes $B$ of $\sub{\mathcal{G}}{\Sigma}$ by putting $B^{(\bullet)} \defAs \{q^{(\bullet)} \st q \in B\}$.

Having defined the vertex set of $\sub{\mathcal{G}}{\Sigma}$, next we describe its edge set. The edges issuing from each place $q$ are exactly all pairs $\langle q,B \rangle$ such that $q\in B \in \sub{\Nodes}{\Sigma}$ (membership edges of $\sub{\mathcal{G}}{\Sigma}$). The remaining edges of $\sub{\mathcal{G}}{\Sigma}$ go from nodes to places (distribution or saturation edges of $\sub{\mathcal{G}}{\Sigma}$). Only places $q$ corresponding to $\otimes$-blocks $q^{(\bullet)}$ of $\Sigma$ (hence called $\otimes$-places) can have incoming edges. Likewise, only nodes $B$ such that $B^{(\bullet)} \in \sub{\Pi}{\otimes}$ (see Definition~\ref{otimes-blocks-upblocks}) can have outgoing edges. Such nodes will be called $\otimes$-nodes.
Specifically, for a $\otimes$-node $B$ and a $\otimes$-place $q$ of $\sub{\mathcal{G}}{\Sigma}$, there is an edge $\langle B,q \rangle$ exactly when 
\[
q^{(\bullet)} \cap \powastot{B^{(\bullet)}} \neq \emptyset,
\]
 namely when there is some ``flow'' of unordered pairs from $\powastot{B^{(\bullet)}}$ to $q^{(\bullet)}$ (through the edge $\langle B,q \rangle$). This is the sense in which a $\otimes$-graph can be considered a kind of flow graph in the realm of set theory. 
Thus, the target map $\sub{\TARGETS}{\Sigma}$ of $\sub{\mathcal{G}}{\Sigma}$ is defined by
\[
\sub{\TARGETS}{\Sigma}(B) \defAs \{q \in \sub{\Places}{\Sigma,\!\otimes} \st q^{(\bullet)} \cap \powastot{B^{(\bullet)}} \neq \emptyset\}, \quad \text{for } B \in \sub{\Nodes}{\Sigma,\!\otimes},
\]
where $\sub{\Places}{\Sigma,\!{\otimes}}$ and $\sub{\Nodes}{\Sigma,\!\otimes}$ denote the collections of the $\otimes$-places and of the $\otimes$-nodes of $\sub{\mathcal{G}}{\Sigma}$, respectively.

Notice that each $\otimes$-node $B$ of $\sub{\mathcal{G}}{\Sigma}$ has some target. Indeed, from $\medcup \sub{\Sigma}{\otimes} = \medcup \powAstot[\sub{\Pi}{\otimes}]$ it follows that $\emptyset \neq \powastot{B^{(\bullet)}} \subseteq \medcup \sub{\Sigma}{\otimes}$, and therefore $\sub{\TARGETS}{\Sigma}(B) \neq \emptyset$.

\subsubsection{Accessible $\otimes$-graphs}
Only \emph{accessible} $\otimes$-graphs are relevant for our decidability purposes. 

\begin{mydef}[Accessible $\otimes$-graphs]\label{def:accessibility}\rm
A place of a $\otimes$-graph $\mathcal{G}=(\Places,\Nodes,\TARGETS)$ is a \textsc{source place} if it has no incoming edges. The remaining places, namely those with incoming edges, are called \textsc{$\otimes$-places}. We denote by $\sub{\Places}{\otimes}$ the collection of the $\otimes$-places of $\mathcal{G}$.

A place of $\mathcal{G}$ is \textsc{accessible} (from the source places of $\mathcal{G}$) if either it is a source place or, recursively, it is the target of some node of $\mathcal{G}$ whose places are \emph{all} accessible from the source places of $\mathcal{G}$. Finally, a $\otimes$-graph is \textsc{accessible} when all its places are accessible.\footnote{Thus, a $\otimes$-graph with no source places is trivially not accessible.}
\end{mydef}

The following result holds.
\begin{mylemma}\label{le:accessibility}
The $\otimes$-graph $\sub{\mathcal{G}}{\Sigma}$ induced by a given partition $\Sigma$ is accessible.
\end{mylemma}
\begin{proof}
Let $\sub{\mathcal{G}}{\Sigma} = (\sub{\Places}{\Sigma},\sub{\Nodes}{\Sigma},\sub{\TARGETS}{\Sigma})$ be the $\otimes$-graph induced by the partition $\Sigma$ via a given bijection $q \mapsto q^{(\bullet)}$ from $\sub{\Places}{\Sigma}$ onto $\Sigma$.

For contradiction, let us assume that $\sub{\mathcal{G}}{\Sigma}$ is not accessible. Among the non-accessible places of $\sub{\mathcal{G}}{\Sigma}$, we select a place $q \in \sub{\Places}{\Sigma}$ whose corresponding block $q^{(\bullet)} \in \subM{\Sigma}$ contains an element $s$ of smallest rank. Plainly, $q^{(\bullet)}$ must be a $\otimes$-block, because otherwise $q$ would be a source place and therefore it would be trivially accessible. Thus, $q^{(\bullet)} \subseteq \medcup\powAstot[\sub{\Pi}{\otimes}]$, so $s \in \medcup\powAstot[\sub{\Pi}{\otimes}]$, where $\sub{\Pi}{\otimes}$ is the collection of the $\otimes$-upblocks of $\Sigma$. Hence, $s = \{s_{1}, s_{2}\} \in \powastot{B^{(\bullet)}}$, for some $\otimes$-node $B = \{q_{1},q_{2}\}$ such that $s_{1} \in q_{1}^{(\bullet)}$ and $s_{2} \in q_{2}^{(\bullet)}$, and therefore $q \in \sub{\TARGETS}{\Sigma}(B)$. Since $q_{1}^{(\bullet)}$ and $q_{2}^{(\bullet)}$ contain elements of rank strictly less than the rank of $s$, the places $q_{1}$ and $q_{2}$ must be accessible. Thus, after all, the place $q$ would be one of the targets of a node whose places are both accessible, and therefore it would be accessible, contradicting our assumption. Hence, $\sub{\mathcal{G}}{\Sigma}$ is  accessible.
\end{proof}

\subsection{Fulfillment by an accessible $\otimes$-graph}\label{Fulfillment by accessible graphs}
Our next task is to figure out which additional properties are enjoyed by the $\otimes$-graph $\sub{\mathcal{G}}{\Sigma}$ induced by a partition $\Sigma$ (via a certain bijection $q \mapsto q^{(\bullet)}$ from $\sub{\Places}{\Sigma}$ onto $\Sigma$) that satisfies a given \BSTuC-conjunction $\myphi$.

Thus, let us assume that $\Sigma$ satisfies a conjunction $\myphi$ via a partition assignment $\mathfrak{I} \colon \Vars{\myphi} \rightarrow \pow{\Sigma}$. Our sought-after properties will pertain the abstraction $\sub{\mathfrak{F}}{\Sigma} \colon \Vars{\myphi} \rightarrow \pow{\sub{\Places}{\Sigma}}$ of the map $\mathfrak{I}$, which is defined by
\[
\sub{\mathfrak{F}}{\Sigma}(x) \defAs \{ q \in \sub{\Places}{\Sigma} \st q^{(\bullet)} \in \mathfrak{I}(x)\}\/, \qquad \text{for $x \in \Vars{\myphi}$.}
\]

For every literal $x=y \star z$ in $\myphi$, with $\star \in \{\cup,\setminus\}$, we have $\mathfrak{I} \models x=y \star z$ and therefore, by Lemma~\ref{partitionAssignmentBoolean}, $\mathfrak{I}(x) = \mathfrak{I}(y) \star \mathfrak{I}(z)$. Hence, $\sub{\mathfrak{F}}{\Sigma}(x) = \sub{\mathfrak{F}}{\Sigma}(y) \star \sub{\mathfrak{F}}{\Sigma}(z)$.

Similarly, for every literal $x \neq y$ in $\myphi$, we have $\mathfrak{I} \models x \neq y$, and therefore (again by Lemma~\ref{partitionAssignmentBoolean}) $\mathfrak{I}(x) \neq \mathfrak{I}(y)$. Thus, we can derive $\sub{\mathfrak{F}}{\Sigma}(x) \neq \sub{\mathfrak{F}}{\Sigma}(y)$.

Finally, let $x= y \otimes z$ be a $\otimes$-literal in $\myphi$, so that 
\begin{equation}\label{otimes-identity}
\medcup \mathfrak{I}(x) = \medcup \mathfrak{I}(y) \otimes \medcup \mathfrak{I}(z)
\end{equation}
holds. 

We preliminarily observe that $\mathfrak{I}(x) \subseteq \sub{\Sigma}{\otimes}$. Indeed, by \eqref{otimes-identity}, $\medcup \mathfrak{I}(x) = \medcup \powAstot[\mathfrak{I}(y) \otimes \mathfrak{I}(z)]$, and therefore $\mathfrak{I}(x)$ is a $\otimes$-subpartition of $\Sigma$.

Let $\upsilon \in \sub{\mathfrak{F}}{\Sigma}(y)$ and $\zeta \in \sub{\mathfrak{F}}{\Sigma}(z)$. Then $\upsilon^{(\bullet)} \subseteq \medcup \mathfrak{I}(y)$ and $\zeta^{(\bullet)} \subseteq \medcup \mathfrak{I}(z)$, and consequently, by \eqref{otimes-identity}, $\upsilon^{(\bullet)} \otimes \zeta^{(\bullet)} \subseteq \medcup \mathfrak{I}(x)$.
Hence,
\begin{align*}
\emptyset &\neq \big\{ q^{(\bullet)} \in \Sigma \st q^{(\bullet)} \cap (\upsilon^{(\bullet)} \otimes \zeta^{(\bullet)}) \neq \emptyset \big\} \subseteq  \mathfrak{I}(x).\\
\intertext{Thus,}
\emptyset &\neq \big\{ q \in \sub{\Places}{\Sigma} \st q^{(\bullet)} \cap (\upsilon^{(\bullet)} \otimes \zeta^{(\bullet)}) \neq \emptyset \big\} \subseteq \sub{\mathfrak{F}}{\Sigma}(x) \cap \sub{\Places}{\otimes},\\
\intertext{so that}
\emptyset &\neq \sub{\TARGETS}{\Sigma}(\{\upsilon, \zeta\}) = \big\{ q \in \sub{\Places}{\Sigma,\!\otimes} \st q^{(\bullet)} \cap (\upsilon^{(\bullet)} \otimes \zeta^{(\bullet)}) \neq \emptyset \big\} \subseteq \sub{\mathfrak{F}}{\Sigma}(x).
\end{align*}

Next, let $q \in \sub{\mathfrak{F}}{\Sigma}(x)$, so that $q^{(\bullet)} \subseteq \medcup \mathfrak{I}(x)$. Let $s \in q^{(\bullet)}$. Hence, by \eqref{otimes-identity}, $s \in \medcup \mathfrak{I}(y) \otimes \medcup \mathfrak{I}(z)$, and therefore $s \in \upsilon^{(\bullet)} \otimes \zeta^{(\bullet)}$, for some $\upsilon,\zeta \in \sub{\Places}{\Sigma}$ such that $\upsilon^{(\bullet)} \in \mathfrak{I}(y)$ and $\zeta^{(\bullet)} \subseteq \mathfrak{I}(z)$. Thus, $\upsilon \in \subM{\mathfrak{F}}(y)$ and $\zeta \in \subM{\mathfrak{F}}(z)$. Since $q^{(\bullet)} \cap (\upsilon^{(\bullet)} \otimes \zeta^{(\bullet)}) \neq \emptyset$ and $q^{(\bullet)} \in \sub{\Sigma}{\otimes}$, we have 
$
q \in \sub{\TARGETS}{\Sigma}(\{\upsilon,\zeta\}) \subseteq \medcup \sub{\TARGETS}{\Sigma}[\sub{\mathfrak{F}}{\Sigma}(y) \otimes \sub{\mathfrak{F}}{\Sigma}(z)].
$
Hence, from the arbitrariness of $q \in \sub{\mathfrak{F}}{\Sigma}(x)$, it follows that
\begin{equation}\label{FSigmax-inclusion}
\sub{\mathfrak{F}}{\Sigma}(x) \subseteq \medcup \sub{\TARGETS}{\Sigma}[\sub{\mathfrak{F}}{\Sigma}(y) \otimes \sub{\mathfrak{F}}{\Sigma}(z)].
\end{equation}

Finally, we prove that also the following identity holds:
\begin{equation}\label{emptyIntersection}
\medcup \sub{\TARGETS}{\Sigma}[\sub{\Nodes}{\Sigma} \setminus (\sub{\mathfrak{F}}{\Sigma}(y) \otimes \sub{\mathfrak{F}}{\Sigma}(z))] \cap \sub{\mathfrak{F}}{\Sigma}(x) = \emptyset.
\end{equation}
For contradiction, let us assume that there exists some $q \in \medcup \sub{\TARGETS}{\Sigma}[\sub{\Nodes}{\Sigma} \setminus (\sub{\mathfrak{F}}{\Sigma}(y) \otimes \sub{\mathfrak{F}}{\Sigma}(z))] \cap \sub{\mathfrak{F}}{\Sigma}(x)$.  Hence, $q \in \sub{\TARGETS}{\Sigma}(A)$, for some $A \in \sub{\Nodes}{\Sigma} \setminus (\sub{\mathfrak{F}}{\Sigma}(y) \otimes \sub{\mathfrak{F}}{\Sigma}(z))$. Since $q^{(\bullet)} \cap \powastot{A^{(\bullet)}} \neq \emptyset$ and $q^{(\bullet)} \subseteq \medcup \mathfrak{I}(x)$, by \eqref{otimes-identity} $\powastot{A^{(\bullet)}} \cap (\medcup \mathfrak{I}(y) \otimes \medcup \mathfrak{I}(z)) \neq \emptyset$, and therefore $A^{(\bullet)} \in \mathfrak{I}(y) \otimes \mathfrak{I}(z)$. The latter membership implies $A \in \sub{\mathfrak{F}}{\Sigma}(y) \otimes \sub{\mathfrak{F}}{\Sigma}(z)$, which is a contradiction, thus proving \eqref{emptyIntersection}.

We can summarize what we have just proved by saying that the accessible $\otimes$-graph induced by a partition $\Sigma$ satisfying a given \BSTuC-conjunction $\myphi$ \emph{fulfills} $\myphi$, according to the following definition.

\begin{mydef}[Fulfillment by an accessible $\otimes$-graph]\label{satAccessible}\rm
An accessible $\otimes$-graph $\mathcal{G}=(\Places,\Nodes,\TARGETS)$  \textsc{fulfills} a given \BSTuC-conjunction $\myphi$ provided that there exists a map $\mathfrak{F} \colon \Vars{\myphi} \rightarrow \pow{\Places}$ (called a \textsc{$\mathcal{G}$-fulfilling map for $\myphi$})  such that the following conditions are satisfied:
\begin{enumerate}[label=(\alph*), ref=\alph*]
\item\label{satAccessibleA} $\mathfrak{F}(x) = \mathfrak{F}(y) \star \mathfrak{F}(z)$, for every conjunct $x=y \star z$ in $\myphi$, where $\star \in \{\cup,\setminus\}$;

\item\label{satAccessibleB} $\mathfrak{F}(x) \neq \mathfrak{F}(y)$, for every conjunct $x \neq y$ in $\myphi$;

\item\label{satAccessibleC} for every conjunct $x=y \otimes z$ in $\myphi$,

\begin{enumerate}[label=(\ref{satAccessibleC}$_{\arabic*}$)]
\item\label{c1} $\emptyset \neq \Targets{ \{\upsilon,\zeta\} } \subseteq \mathfrak{F}(x)$, for all $\upsilon \in \mathfrak{F}(y)$ and $\zeta \in \mathfrak{F}(z)$;

\item\label{c2} $\mathfrak{F}(x) \subseteq \bigcup\TARGETS[\mathfrak{F}(y) \otimes \mathfrak{F}(z)]$;

\item\label{c3} $\bigcup\TARGETS[\Nodes \setminus (\mathfrak{F}(y) \otimes \mathfrak{F}(z))] \cap \mathfrak{F}(x) = \emptyset$.
\end{enumerate}
\end{enumerate}

\end{mydef}

Thus, the results just proved can be stated as follows.

\begin{mylemma}\label{otimesGraphFulfillsMyphi}
The accessible $\otimes$-graph induced by a partition satisfying a given \BSTuC-conjunction $\myphi$ fulfills $\myphi$.
\end{mylemma}

As an immediate consequence, we have:

\begin{mycorollary}\label{cor-otimesGraphFulfillsMyphi}
A satisfiable \BSTuC-conjunction with $n$ variables is fulfilled by an accessible $\otimes$-graph of size (at most) $2^{n}-1$.
\end{mycorollary}
\begin{proof}
Let $\myphi$ be a satisfiable \BSTuC-conjunction with $n$ variables.
As stated in Lemma~\ref{wasA}, $\myphi$ is satisfied by a partition $\Sigma$ with exactly $2^{n}-1$ blocks. Thus, the $\otimes$-graph $\sub{\mathcal{G}}{\Sigma}$ induced by $\Sigma$ has size $2^{n}-1$ and, by Lemmas~\ref{le:accessibility} and \ref{otimesGraphFulfillsMyphi}, it is accessible and fulfills $\myphi$.
\end{proof}

\subsection{Construction process}
Lemma~\ref{otimesGraphFulfillsMyphi} can be reversed, thus yielding a decision procedure for \BSTuC-conjunctions.

\begin{mylemma}\label{rev-otimesGraphFulfillsMyphi}
If a \BSTuC-conjunction is fulfilled by an accessible $\otimes$-graph, then it is satisfiable.
\end{mylemma}
\begin{proof}
Let $\mathcal{G}=(\Places,\Nodes,\TARGETS)$ be an accessible $\otimes$-graph, and let us assume that $\mathcal{G}$ fulfills a given \BSTuC-conjunction $\myphi$ via the map $\mathfrak{F} \colon \Vars{\myphi} \rightarrow \pow{\Places}$. 

To each place $q \in \Places$, we associate a set $q^{(\bullet)}$, initially empty. Then, by suitably exploiting the $\otimes$-graph $\mathcal{G}$ as a kind of flow graph, we shall show that the sets $q^{(\bullet)}$'s can be monotonically extended by a (possibly infinite) \emph{construction process} (comprising a finite \emph{initialization phase} and a subsequent (possibly infinite) \emph{stabilization phase}) in such a way that the following properties hold:
\begin{enumerate}[label=(P\arabic*)]
\item\label{P1} After each step, the sets $q^{(\bullet)}$'s are pairwise disjoint.

\item\label{P2} 
At the end of the initialization phase all the $q^{(\bullet)}$'s are nonempty (and pairwise disjoint). Thus, after each step in the subsequent stabilization phase,  the sets $q^{(\bullet)}$'s, with  $q \in \Places$, form a partition equipollent with $\Places$.

\item\label{P3} 
After each step in the stabilization phase, the inclusion 
\[
q^{(\bullet)} \subseteq \medcup \big\{ \powastot{A^{(\bullet)}} \st A \in \TARGETS^{-1}(q)\big\}
\]
holds, for each $\otimes$-place $q \in \sub{\Places}{\otimes}$, where we are using the notation $B^{(\bullet)} \defAs \{p^{(\bullet)} \st p \in B\}$ for $B \in \Nodes$.
\end{enumerate}

\begin{enumerate}[resume*]
\item\label{P4} At the end of the construction process, we have 
\[
\powastot{A^{(\bullet)}} \subseteq \bigcup \{q^{(\bullet)} \st q \in \Targets{A}\},
\]
for each $\otimes$-node $A \in \sub{\Nodes}{\otimes}$ (namely for each node $A \in \Nodes$ such that $\Targets{A} \neq \emptyset$).\footnote{Should the construction process involve denumerably many steps, the final values of the $q^{(\bullet)}$'s are to be intended as limit of the sequences of their values after each step in the stabilization phase.}
\end{enumerate}

Subsequently, we shall prove that the properties \ref{P1}--\ref{P4} together with the conditions \eqref{satAccessibleA}--\eqref{satAccessibleC} of Definition~\ref{satAccessible}, characterizing the fulfilling $\otimes$-graph $\mathcal{G}$, allow one to show that the partition $\{q^{(\bullet)} \st q \in \Places\}$ resulting from the above construction process satisfies our conjunction $\myphi$.

\bigskip

The initialization and stabilization phases of our construction process consist of the following steps.
\paragraph{Initialization phase:}
\begin{enumerate}[label=(I$_{\arabic*}$)]
\item\label{I1} To begin with, let $\{\overline q \st q \in \Places \setminus \sub{\Places}{\otimes}\}$ be any partition equipollent to the set $\Places \setminus \sub{\Places}{\otimes}$ of the source places of $\mathcal{G}$, where each block $\overline q$, for $q \in \Places \setminus \sub{\Places}{\otimes}$, is a hereditarily finite set of cardinality (at least) $\max(2|\sub{\Places}{\otimes}|,1)$ and whose members all have cardinality \emph{strictly} greater than 2,\footnote{For the present case concerning the ordinary satisfiability problem, we could have allowed that the $\overline q$'s were all infinite sets, rather than hereditarily finite sets. However, we chose to enforce hereditarily finiteness of the $\overline q$'s even in the current case in order that the initialization phase would coincide with that for the hereditarily finite satisfiability case to be addressed in the next section.} and put
\[
q^{(\bullet)} \defAs \begin{cases}
\overline q & \text{if } q \in \Places \setminus \sub{\Places}{\otimes}\\
\emptyset & \text{if } q \in \sub{\Places}{\otimes}.
\end{cases}
\]
\end{enumerate}

We say that a place $q \in \Places$ has already been \emph{initialized} when $q^{(\bullet)} \neq \emptyset$. Likewise, a $\otimes$-node $A \in \sub{\Nodes}{\otimes}$ has been \emph{initialized} when its places have all been initialized. During the initialization phase, an initialized $\otimes$-node $A \in \sub{\Nodes}{\otimes}$ is said to be \textsc{ready} if it has some target that has not been yet initialized. 

\begin{enumerate}[resume*]
\item\label{I2} While there are places in $\Places$ not yet initialized, pick any ready node $A \in \Nodes$ and distribute evenly all the members of $\powastot{A^{(\bullet)}}$ among all of its targets.
\end{enumerate}

The accessibility of $\mathcal{G}$ guarantees that the while-loop \ref{I2} terminates in a finite number of iterations.

At the end of the initialization phase all the $q^{(\bullet)}$'s are nonempty, so  property \ref{P2} holds. Indeed, if there were no $\otimes$-places, then all places would be initialized just after step \ref{I1}, and so all the $q^{(\bullet)}$'s would be (nonempty) singletons. On the other hand, if $|\sub{\Places}{\otimes}|>0$, then at the end of the while-loop \ref{I2} we shall have $|q^{(\bullet)}| \geq 2|\sub{\Places}{\otimes}|$, for each $q \in \Places$. This follows just from the initialization step \ref{I1}, for all source places $q \in \Places \setminus \sub{\Places}{\otimes}$. Otherwise, by induction, we have $|q^{(\bullet)}| \geq 2|\sub{\Places}{\otimes}|$, for every $q$ in a ready node $A \in \Nodes$, and therefore 
\[
|\powastot{A^{(\bullet)}}| \geq {2|\sub{\Places}{\otimes}| \choose 2} + 2|\sub{\Places}{\otimes}| = |\sub{\Places}{\otimes}| \cdot (2|\sub{\Places}{\otimes}| + 1).
\]
Hence, each of the $|\Targets{A}| \leq |\sub{\Places}{\otimes}|$ sets $t^{(\bullet)}$, for $t \in \Targets{A}$, will receive at least $2|\sub{\Places}{\otimes}| + 1$ elements by the distribution step relative to the node $A$.

Concerning property \ref{P1}, we observe that at each distribution step, only elements of cardinality 1 or 2 are added to the sets $q^{(\bullet)}$'s. Therefore, the disjointness of the sets $q_{1}^{(\bullet)}$ and $q_{2}^{(\bullet)}$, for any two distinct places $q_{1},q_{2} \in \Places$ such that at least one of them is a source place, will be guaranteed. Indeed, if both $q_{1}$ and $q_{2}$ are source places, then $q_{1}^{(\bullet)} \cap q_{2}^{(\bullet)} = \overline q_{1} \cap \overline q_{2} = \emptyset$. On the other hand, if only one of them is a source node, say $q_{1}$, then since $q_{1}^{(\bullet)} = \overline q_{1}$ contains only members of cardinality strictly greater than 2 whereas, by step \ref{I2}, all the members of $q_{2}^{(\bullet)}$ have cardinality less than or equal to 2, it follows that even in this case we have $q_{1}^{(\bullet)} \cap q_{2}^{(\bullet)} = \emptyset$. Finally, for any two distinct places $q_{1}$ and $q_{2}$, none of which is a source node, we observe that if they have been initialized by a distribution step applied to the same node $A \in \Nodes$, we trivially have $q_{1}^{(\bullet)} \cap q_{2}^{(\bullet)} = \emptyset$. On the other hand, if $q_{i}$ is initialized by distributing over a $\otimes$-node $A_{i}$, with $i=1,2$, where $A_{1} \neq A_{2}$, by induction on the number of the number of distribution steps it can easily be shown that $A_{1}^{(\bullet)} \neq A_{2}^{(\bullet)}$, and therefore $q_{1}^{(\bullet)} \cap q_{2}^{(\bullet)} \subseteq \powastot{A_{1}^{(\bullet)}} \cap \powastot{A_{2}^{(\bullet)}} = \emptyset$.

\paragraph{Stabilization phase:} During the stabilization phase, a $\otimes$-node $A \in \sub{\Nodes}{\otimes}$ is \textsc{ripe} if 
\[
\powastot{A^{(\bullet)}} \setminus \medcup \big\{q^{(\bullet)} \st q \in \Targets{A}\big\} \neq \emptyset.
\]

We execute the following (possibly infinite) loop:
\begin{enumerate}[label=(S$_{\arabic*}$)]
\item\label{while-loop} While there are ripe $\otimes$-nodes, pick any of them, say $A \in \Nodes$, and distribute all the members of $\powastot{A^{(\bullet)}} \setminus \medcup \big\{q^{(\bullet)} \st q \in \Targets{A}\big\}$ (namely the members of $\powastot{A^{(\bullet)}}$ that have not been distributed yet) among its targets howsoever.
\end{enumerate}

The fairness condition that one must comply with is the following: 
\begin{quote}
once a $\otimes$-node becomes ripe during the stabilization phase, it must be picked for distribution within a finite number of iterations of the while-loop \ref{while-loop}.
\end{quote}
A possible way to enforce such condition consists, for instance, in maintaining all ripe $\otimes$-nodes in a queue $\mathcal{Q}$, picking always the $\otimes$-node to be used in a distribution step from the front of $\mathcal{Q}$ and adding the $\otimes$-nodes that have just become ripe to the back of $\mathcal{Q}$, provided that they are not already in $\mathcal{Q}$.

By induction on $n \in \Nats$, it is not hard to show that properties \ref{P1} and \ref{P3} will hold just after the $n$-th iteration of the while-loop \ref{while-loop} of the stabilization phase, and that property \ref{P4} will hold at the end of the stabilization phase, in case of termination.

Instead, when the stabilization phase runs for denumerably many steps, the final partition $\Places^{(\bullet)}$ is to be intended as the limit of the partial partitions constructed after each step of the stabilization phase. Specifically, for each place $q \in \Places$, we let $q^{(i)}$ be the value of $q^{(\bullet)}$ just after the $i$-th iteration of \ref{while-loop}. Plainly, we have 
\begin{equation}\label{monotonicity}
q^{(i)} \subseteq q^{(i+1)}, \quad \text{for } i \in \Nats. 
\end{equation}
Then we put
\begin{equation}\label{q-bullet-limit}
q^{(\bullet)} \defAs \bigcup_{i \in \Nats} q^{(i)}, \quad \text{for } q \in \Places
\end{equation}
(notation overloading should not be a problem).

By way of illustration, we prove that property \ref{P4} holds for the partition $\Places^{(\bullet)} = \big\{q^{(\bullet)} \st q \in \Places \big\}$, when the $q^{(\bullet)}$'s are defined by \eqref{q-bullet-limit}. To this purpose, let $A \in \Nodes$ be such that $\Targets{A} \neq \emptyset$, and assume for contradiction that 
\[
\powastot{A^{(\bullet)}} \not\subseteq \medcup \big\{q^{(\bullet)} \st q \in \Targets{A}\big\}.
\]
Let $s$ be any element in $\powastot{A^{(\bullet)}} \setminus \medcup \big\{q^{(\bullet)} \st q \in \Targets{A}\big\}$, and let $i \in \Nats$ be the smallest index such that $s \in \powastot{A^{(i)}}$, where $A^{(i)} \defAs \{q^{(i)} \st q \in A \}$. Since $s \in \powastot{A^{(i)}} \setminus \medcup \big\{q^{(i)} \st q \in \Targets{A}\big\}$, the node $A$ must have been ripe just after the $i$-th iteration of \ref{while-loop}. Therefore, by the fairness condition, the node $A$ will be picked for distribution in a finite number of steps, say $k$, after the $i$-th step, so that we have
\begin{align*}
\powastot{A^{(i)}} &\subseteq \powastot{A^{(i+k)}} && \text{(by \eqref{monotonicity})}\\
&\subseteq \medcup \big\{q^{(i+k+1)} \st q \in \Targets{A}\big\}\\
&\subseteq \medcup \big\{q^{(\bullet)} \st q \in \Targets{A}\big\},
\end{align*}
and therefore $s \in \medcup \big\{q^{(\bullet)} \st q \in \Targets{A}\big\}$, which is a contradition. Thus, property \ref{P4} holds also when the construction process takes a denumerable number of steps.

\bigskip

Next, we show that the final partition $\Places^{(\bullet)} = \{q^{(\bullet)} \st q \in \Places \}$ satisfies $\myphi$. In particular, we prove that the partition assignment $\mathfrak{I} \colon \Vars{\myphi} \rightarrow \pow{\Places^{(\bullet)}}$ defined by
\[
\mathfrak{I}(x) \defAs \{q^{(\bullet)} \st q \in \mathfrak{F}(x) \} \text{, \quad for } x \in \Vars{\myphi},
\]
satisfies $\myphi$, where we recall that $\mathfrak{F}$ is the $\mathcal{G}$-fulfilling map for $\myphi$.

Since $\mathfrak{F}$ is a $\mathcal{G}$-fulfilling map for $\myphi$, then
\begin{enumerate}[label=-]
\item for every literal $x=y \star z$ in $\myphi$, with $\star \in \{\cup,\setminus\}$, we have $\mathfrak{F}(x) = \mathfrak{F}(y) \star \mathfrak{F}(z)$, so that $\mathfrak{I}(x) = \mathfrak{I}(y) \star \mathfrak{I}(z)$ holds; and

\item for every literal $x \neq y$ in $\myphi$, we have $\mathfrak{F}(x) \neq \mathfrak{F}(y)$, so that $\mathfrak{I}(x) \neq \mathfrak{I}(y)$ holds.
\end{enumerate}
Thus, by Lemma~\ref{partitionAssignmentBoolean}, the partition assignment $\mathfrak{I}$ satisfies all Boolean literals in $\myphi$ of types
\[
x=y \cup z, \quad x=y \setminus z, \quad x \neq y.
\]

Next, let $x=y \otimes z$ be a conjunt of $\myphi$. We prove separately that the following inclusions hold:
\begin{gather}
\medcup \mathfrak{I}(x) \subseteq \medcup \mathfrak{I}(y) \otimes \medcup \mathfrak{I}(z) \label{firstInclusion}\\
\medcup \mathfrak{I}(y) \otimes \medcup \mathfrak{I}(z) \subseteq \medcup \mathfrak{I}(x).\label{secondInclusion}
\end{gather}

Concerning \eqref{firstInclusion}, let $q^{(\bullet)} \subseteq \bigcup\mathfrak{I}(x)$. Then $q^{(\bullet)} \in \mathfrak{I}(x)$, so that $q \in \mathfrak{F}(x)$. By \ref{c2}, $q$ cannot be a source place. Hence, by \ref{P3}, we have:
\[
q^{(\bullet)} \subseteq \medcup \big\{ \powastot{A^{(\bullet)}} \st A \in \TARGETS^{-1}(q)\big\}.
\]
Next we show that 
\begin{equation}\label{threeStars}
\TARGETS^{-1}(q) \subseteq \mathfrak{F}(y) \otimes \mathfrak{F}(z).
\end{equation}
Let $A \in \TARGETS^{-1}(q)$ (so that $q \in \Targets{A}$), and for contradiction assume that $A \notin \mathfrak{F}(y) \otimes \mathfrak{F}(z)$. Then, by \ref{c3}, we have $\Targets{A} \cap \mathfrak{F}(x) = \emptyset$, contradicting $q \in \Targets{A} \cap \mathfrak{F}(x)$. Thus, $A \in \mathfrak{F}(y) \otimes \mathfrak{F}(z)$, proving \eqref{threeStars}. Hence, we have:
\begin{align*}
q^{(\bullet)} &\subseteq \medcup \big\{ \powastot{A^{(\bullet)}} \st A \in \TARGETS^{-1}(q)\big\}\\
&\subseteq \medcup \big\{ \powastot{A^{(\bullet)}} \st A \in \mathfrak{F}(y) \otimes \mathfrak{F}(z)\big\}\\
&= \medcup \big\{ \powastot{A^{(\bullet)}} \st A^{(\bullet)} \in \mathfrak{I}(y) \otimes \mathfrak{I}(z)\big\}\\
&= \medcup \mathfrak{I}(y) \otimes \medcup \mathfrak{I}(z) && \text{(by Lemma~\ref{wasD}),}
\end{align*}
and therefore the inclusion \eqref{firstInclusion} holds.

Concerning the inclusion \eqref{secondInclusion}, let $s \in \medcup \mathfrak{I}(y) \otimes \medcup \mathfrak{I}(z)$. Hence, $s \in q_{1}^{(\bullet)} \otimes q_{2}^{(\bullet)} = \powastot{\{q_{1}^{(\bullet)},q_{2}^{(\bullet)}\}}$, for some $q_{1} \in \mathfrak{F}(y)$ and $q_{2} \in \mathfrak{F}(z)$. From \ref{c1}, we have $\emptyset \neq \Targets{\{q_{1},q_{2}\}} \subseteq \mathfrak{F}(x)$. Thus, by \ref{P4},
\begin{align*}
\powastot{\{q_{1}^{(\bullet)},q_{2}^{(\bullet)}\}} &\subseteq \medcup \{ q^{(\bullet)} \st q \in \Targets{\{q_{1},q_{2}\}} \}\\
&\subseteq \medcup \{ q^{(\bullet)} \st q \in  \mathfrak{F}(x) \}\\
&= \medcup \mathfrak{I}(x),
\end{align*}
and therefore $s \in \medcup \mathfrak{I}(x)$, proving \eqref{secondInclusion} by the arbitrariness of $s \in \medcup \mathfrak{I}(y) \otimes \medcup \mathfrak{I}(z)$. 

Hence, the partition assignment $\mathfrak{I}$ satisfies also all the literals in $\myphi$ of the form $x=y \otimes z$, and in turn the final partition $\Places^{(\bullet)}$ satisfies the conjunction $\myphi$.
\end{proof}

By combining Lemmas~\ref{otimesGraphFulfillsMyphi} and \ref{rev-otimesGraphFulfillsMyphi} and Corollary~\ref{cor-otimesGraphFulfillsMyphi}, we obtain:
\begin{mytheorem}
A \BSTuC-conjunction with $n$ variables is satisfiable if and only if it is fulfilled by an accessible $\otimes$-graph of size (at most) $2^{n}-1$.
\end{mytheorem}

The preceding theorem is at the base of the following trivial decision procedure for \BSTuC:

\begin{quote}{\small
\begin{tabbing}
xx \= \= xx \= xx \= xx \= xx \= xx \= xx \= xx \= xx \kill
\hspace{-10pt}\textbf{procedure} \textsf{\BSTuC-satisfiability-test}$(\myphi)$;\\
\> \ninstrb \> $n \defAs |\Vars{\myphi}|$;\\
\> \ninstrb \> \textbf{for} each $\otimes$-graph $\mathcal{G}$ with $2^{n}-1$ places \textbf{do}\\
\> \ninstrb \> \> \textbf{if} $\mathcal{G}$ is accessible and fulfills $\myphi$ \textbf{then}\\
\> \ninstrb \> \> \> \textbf{return} ``$\myphi$ is satisfiable'';\\
\> \ninstrb \> \textbf{return} ``$\myphi$ is unsatisfiable'';\\
\hspace{-10pt}\textbf{end procedure};
\end{tabbing}
}
\end{quote}

Concerning the complexity of the above procedure, we observe that, given a \BSTuC-conjunction $\myphi$ with $n$ distinct variables, we have:
\begin{enumerate}[label=-]
\item the size of $\myphi$ is $\mathcal{O}(n^{3})$ (w.l.o.g., we are assuming that literal repetitions are not admitted);%

\item the size of a $\otimes$-graph with $2^{n}-1$ places is $\mathcal{O}(8^{n})$;

\item the size of any candidate fulfilling map over a set of $n$ variables is $\mathcal{O}(n 2^{n})$ and the time needed to check whether it is actually a $\mathcal{G}$-fulfilling map for $\myphi$, for a given $\otimes$-graph $\mathcal{G}$ with $2^{n}-1$ places, is $\mathcal{O}(8^{n})$.
\end{enumerate}

Hence, for a \BSTuC-conjunction $\myphi$ with $n$ distinct variables the procedure \textsf{\BSTuC-satisfiability-test} has a nondeterministic $\mathcal{O}(8^{n})$-time complexity. Hence, we have:

\begin{mytheorem}
The satisfiability problem for \BSTuC-conjunctions belongs to the complexity class \emph{\textsf{NEXPTIME}}.
\end{mytheorem}

The above result can be easily generalized to \BSTuC-formulae that are not necessarily conjunctions.

\begin{mytheorem}\label{th:BSTuC-formulae is in NEXPTIME}
The satisfiability problem for \BSTuC-formulae belongs to the complexity class \emph{\textsf{NEXPTIME}}.
\end{mytheorem}

\smallskip

There are satisfiable \BSTuC-formulae that admit only infinite models. This is the case, for instance, for the following conjunction $\myphi$
\[
x \neq x \setminus x \:\wedge\: x \otimes x \subseteq x,
\]
which is satisfied by the assignment $M$ such that $Mx = \HF$. In addition, for every model $M$ of $\myphi$, the first conjunct $x \neq x \setminus x$ forces $Mx$ to be  nonempty, while the second conjunct $x \otimes x \subseteq x$ forces $Mx$ to be infinite. Indeed, whenever a set $s$ belongs to $Mx$, it must also be the case that its singleton $\{s\}$ belongs to $Mx$ as well. Thus, iteratively, the infinitely many sets 
\[
\{s\},~ \{\{s\}\},~ \{\{\{s\}\}\},~\ldots
\]
must all belong to $s$, proving that $Mx$ must be infinite.

It is therefore important to investigate the finite s.p.\ for \BSTuC, which we do in the next section.

\section{The finite and the hereditarily finite satisfiability problems for \BSTuC}\label{se:finite satisfiability}

Let $\myphi$ be a finitely satisfiable \BSTuC-conjunction, and let now $\Sigma$ be a partition with \emph{finite} domain $\medcup \Sigma$ that satisfies  $\myphi$ via some partition assignment $\mathfrak{I} \colon \Vars{\myphi} \rightarrow \pow{\Sigma}$. Also, let $\sub{\mathcal{G}}{\Sigma} = (\sub{\Places}{\Sigma},\sub{\Nodes}{\Sigma},\sub{\TARGETS}{\Sigma})$ be the $\otimes$-graph induced by  $\Sigma$ via a given bijection $q \mapsto q^{(\bullet)}$. As argued just before Lemma~\ref{le:accessibility}, the graph $\sub{\mathcal{G}}{\Sigma}$ is $\otimes$-accessible and fulfills $\myphi$ via the map $\sub{\mathfrak{F}}{\Sigma} \colon \Vars{\myphi} \rightarrow \pow{\sub{\Places}{\Sigma}}$ induced by $\mathfrak{I}$ and defined by 
\[
\sub{\mathfrak{F}}{\Sigma}(x) \defAs \{ q \in \sub{\Places}{\Sigma} \st q^{(\bullet)} \in \mathfrak{I}(x)\}\/, \qquad \text{for $x \in \Vars{\myphi}$}
\]
(so, $\sub{\mathfrak{F}}{\Sigma}$ is a $\sub{\mathcal{G}}{\Sigma}$-fulfilling map for $\myphi$). 

We shall see that the finiteness of $\medcup \Sigma$ yields a weak kind of acyclicity for the induced $\otimes$-graph $\sub{\mathcal{G}}{\Sigma}$, which is expressed in terms of a restricted form of topological order.

\begin{mydef}
\rm
A \textsc{topological $\otimes$-order} of a $\otimes$-graph $\mathcal{G}=(\Places,\Nodes,\TARGETS)$ is any total order $\prec$ over its set of places $\Places$ such that
\begin{equation}\label{max-prec-max}
\max_{\prec} A ~\prec~ \max_{\prec} \TARGETS(A), 
\end{equation}
for every $\otimes$-node $A$ of $\mathcal{G}$.

We write $\mathcal{G}=(\Places,\Nodes,\TARGETS,\prec)$ for a $\otimes$-graph $(\Places,\Nodes,\TARGETS)$ endowed with a topological $\otimes$-order $\prec$, and we refer to it as a \textsc{(topologically) $\otimes$-ordered graph}.
\end{mydef}

Notice that a $\otimes$-ordered graph need not be acyclic. On the other hand, any acyclic $\otimes$-graph admits a topological order of its vertices and therefore a topological $\otimes$-order, as can be easily checked. In this sense, topological $\otimes$-orders are less demanding than ordinary topological orders.

Later we shall also see that, together with fulfillability and accessibility, the existence of a topological $\otimes$-order is sufficient for a \BSTuC-conjunction to be hereditarily finitely satisfiable, thereby proving that the finite and the hereditarily finite satisfiability problems for \BSTuC are equivalent.

To start with, we show, as announced before, that the induce $\otimes$-graph $\sub{\mathcal{G}}{\Sigma}$ admits a topological $\otimes$-order.
Thus, let $\sub{\prec}{\Sigma}$ be any total order over $\sub{\Places}{\Sigma}$ that refines the partial order induced by the rank function, namely such that 
\[
\rk p^{(\bullet)} < \rk q^{(\bullet)} \quad \Longrightarrow \quad p \sub{\prec}{\Sigma} q, \quad \text{for } p,q \in \sub{\Places}{\Sigma}.
\]
We prove that \eqref{max-prec-max} holds for $\sub{\prec}{\Sigma}$, namely $\sub{\prec}{\Sigma}$ is a topological $\otimes$-order of $\sub{\mathcal{G}}{\Sigma}$. So, let $A$ be any $\otimes$-node of $\sub{\mathcal{G}}{\Sigma}$. For each $q \in A$, we select an $s_{q} \in q^{(\bullet)}$ of maximal rank, which exists since $q^{(\bullet)}$ is finite, and put $s_{A} \defAs \{s_{q} \st q \in A\}$. Let $q_{A}$ be the target if $A$ such that $s_{A} \in q_{A}^{(\bullet)}$ (plainly, such a target exists, since $s_{A} \in \powastot{A^{(\bullet)}} \subseteq \medcup \sub{\Sigma}{\otimes}$). Hence, for each $q \in A$ we have
\[
\rk q^{(\bullet)} < \rk s_{A}  < \rk q_{A}^{(\bullet)},
\]
so that $q \sub{\prec}{\Sigma} q_{A}$ holds. But then
\[
\max_{\sub{\prec}{\Sigma}} A ~\sub{\prec}{\Sigma}~ q_{A} ~\subAlt{\preceqq}{\Sigma}~\max_{~\sub{\prec}{\Sigma}~} \TARGETS(A),
\]
proving that $\sub{\prec}{\Sigma}$ is a topological $\otimes$-order of $\sub{\mathcal{G}}{\Sigma}$.

Summing up, we have proved that:
\begin{mylemma}\label{lemma14}
A finitely satisfiable \BSTuC-conjunction is fulfilled by an accessible ordered $\otimes$-graph.
\end{mylemma}

Next, we prove that if a \BSTuC-conjunction $\myphi$ is fulfilled by an accessible ordered $\otimes$-graph, then it is satisfiable by a \emph{hereditarily} finite model.

Thus, let $\mathcal{G} = (\Places, \Nodes, \TARGETS, \prec)$ be an accessible ordered $\otimes$-graph that fulfills $\myphi$ via a map $\mathfrak{F} \colon \Vars{\myphi} \rightarrow \pow{\Places}$, and let $\preceq$ be the \emph{total preorder} induced by $\prec$ over $\Nodes$, defined by
\[
A \preceq B ~~\xleftrightarrow{\hbox{\smaller[5]{$~\mathit{Def}$~}}}~~ \max_{\prec} A ~\preceqq~ \max_{\prec} B,
\]
for all $A,B \in \Nodes$.

Much the same construction process described at depth in the proof of Lemma~\ref{rev-otimesGraphFulfillsMyphi} concerning the ordinary satisfiability problem for \BSTuC will allow us to build a hereditarily finite model for $\myphi$.

Specifically, the initialization phase of our new construction process coincides with that of the old construction process, and therefore consists in the steps \ref{I1} and \ref{I2} seen previously. Instead, the old stabilization loop \ref{while-loop} is replaced by the following one:
\begin{enumerate}[label=(S$_{\arabic*}'$)]
\item\label{while-loop-prime} While there are ripe $\otimes$-nodes, pick any $\preceq$-minimal ripe $\otimes$-node, say $A \in \sub{\Nodes}{\otimes}$, and assign all the members of $\powastot{A^{(\bullet)}} \setminus \medcup \big\{q^{(\bullet)} \st q \in \Targets{A}\big\}$ to the block $q_{A}^{(\bullet)}$ such that $q_{A} = \max_{\prec} \Targets{A}$, namely execute the assignment
\[
q_{A}^{(\bullet)} \coloneqq q_{A}^{(\bullet)} \cup \powastot{A^{(\bullet)}} \setminus \medcup \big\{q^{(\bullet)} \st q \in \Targets{A}\big\}.
\]
(As before, during the stabilization phase a $\otimes$-node $A$ is \emph{ripe} if the set $\powastot{A^{(\bullet)}} \setminus \medcup \big\{q^{(\bullet)} \st q \in \Targets{A}\big\}$ is nonempty.)
\end{enumerate}

We prove that the while-loop \ref{while-loop-prime} can be executed at most $|\sub{\Nodes}{\otimes}|$ times. Thus, let 
\begin{equation}\label{eq:sequence-picked}
A_{1},~A_{2},\ldots,~A_{k},\ldots
\end{equation}
be the sequence of the $\otimes$-nodes picked for distribution during the execution of the loop \ref{while-loop-prime}. It is enough to show that the $\otimes$-nodes in the sequence \eqref{eq:sequence-picked} are pairwise distinct. To this end, we first prove that we have
\begin{equation}\label{eq:preceq-picked}
A_{1} \preceq A_{2}\preceq \ldots \preceq A_{k}\preceq \ldots
\end{equation}
For contradiction, let us assume that \eqref{eq:preceq-picked} does not hold, and let $\ell \in \Nats$ be the least index such that we have
\begin{equation}\label{eq:preceq-contradicts}
A_{\ell} \not\preceq A_{\ell+1},
\end{equation}
so that $A_{\ell+1} \preceq A_{\ell}$ must hold, since the preorder $\preceq$ is total. Plainly, at the $\ell$-th iteration of \ref{while-loop-prime}, the node $A_{\ell+1}$ cannot be ripe, as otherwise it would have been chosen at step $\ell$ in place of $A_{\ell}$. So, the target $q_{A_{\ell}} = \max_{\prec} \TARGETS(A_{\ell})$ of $A_{\ell}$ must belong to $A_{\ell+1}$, and therefore
\[
\max_{\prec} A_{\ell} ~\prec~ \max_{\prec} \TARGETS(A_{\ell}) = q_{A_{\ell}} ~\preceqq~ \max_{\prec} A_{\ell +1}
\]
must hold, yielding $A_{\ell} \preceq A_{\ell +1}$ which contradicts \eqref{eq:preceq-contradicts}.

In what follows, for any node $A \in \Nodes$ we shall denote by $A^{(i)}$ the value of the set $A^{(\bullet)}$ (associated with $A$) just before the $i$-th iteration of the loop \ref{while-loop-prime}.

We are now ready to prove that the nodes in the sequence \eqref{eq:sequence-picked} are pairwise distinct.  For contradiction, if $A_{i} = A_{j}$, with $i < j$, then $A_{i}^{(i)} \neq A_{i}^{(j)}$, so that at least one place $q$ in $A_{i}$ must be the $\prec$-maximum target of some $\otimes$-node, say $A_{t}$ (with $i \leq t \leq j-1$), in the sequence $A_{i},\ldots,A_{j-1}$. But then we would have:
\[
\max_{\prec} A_{i} ~\preceqq~ \max_{\prec} A_{t} ~\prec~ q_{A_{t}} ~\preceqq~ \max_{\prec} A_{i},
\]
which is a contradiction. Therefore, the while-loop \ref{while-loop-prime} must terminate in at most $|\sub{\Nodes}{\otimes}|$ iterations.

Thus, at the end of the construction process under consideration, all the sets $q^{(\bullet)}$, for $q \in \Places$, are plainly hereditarily finite, and since the loop \ref{while-loop-prime} is a particular instance (which is guaranteed to terminate) of the loop \ref{while-loop}, then the partition $\{q^{(\bullet)} \st q \in \Places\}$ resulting from the above construction process satisfies our conjunction $\myphi$, just as argued in the proof of Lemma~\ref{rev-otimesGraphFulfillsMyphi}.

In conclusion, we have:
\begin{mylemma}\label{lemma15}
A \BSTuC-conjunction fulfilled by an accessible ordered $\otimes$-graph is satisfiable by a \emph{hereditarily} finite model.
\end{mylemma}

From Lemmas~\ref{lemma14} and \ref{lemma15} and Corollary~\ref{cor-otimesGraphFulfillsMyphi}, we deduce:
\begin{mytheorem}
The finite and the hereditarily finite satisfiability problems for \BSTuC-conjunctions are equivalent.

In addition, any \BSTuC-conjunction with $n$ variables is (hereditarily) finitely satisfiable if and only if it is fulfilled by an accessible ordered $\otimes$-graph of size (at most) $2^{n}-1$.
\end{mytheorem}

The preceding theorem justifies the following trivial decision procedure for the (hereditarily) finite satisfiability problem for \BSTuC:
\setcounter{instrb}{0}
\begin{quote}{\small
\begin{tabbing}
xx \= \= xx \= xx \= xx \= xx \= xx \= xx \= xx \= xx \kill
\hspace{-10pt}\textbf{procedure} \textsf{\BSTuC-finite-satisfiability-test}$(\myphi)$;\\
\> \ninstrb \> $n \defAs |\Vars{\myphi}|$;\\
\> \ninstrb \> \textbf{for} each $\otimes$-graph $\mathcal{G}$ with $2^{n}-1$ places \textbf{do}\\
\> \ninstrb \> \> \textbf{if} $\mathcal{G}$ is $\otimes$-ordered, accessible and fulfills $\myphi$ \textbf{then}\\
\> \ninstrb \> \> \> \textbf{return} ``$\myphi$ is (hereditarily) finitely satisfiable'';\\
\> \ninstrb \> \textbf{return} ``$\myphi$ is satisfiable by any (hereditarily) finite model'';\\
\hspace{-10pt}\textbf{end procedure};
\end{tabbing}
}
\end{quote}

Much as in the previous section, we can deduce that:
\begin{mytheorem}\label{th:BSTuC-finite sat. is in NEXPTIME}
The (hereditarily) finite satisfiability problem for \BSTuC-formulae belongs to the complexity class \emph{\textsf{NEXPTIME}}.
\end{mytheorem}

We say that a \BSTuC-formula \emph{forces infinite models} if it is satisfiable but not finitely satisfiable.
 
On account of Theorems~\ref{th:BSTuC-formulae is in NEXPTIME} and~\ref{th:BSTuC-finite sat. is in NEXPTIME}, we may finally infer the following immediate result:
\begin{mycorollary}
The problem of deciding whether a \BSTuC-formula forces infinite models belongs to the complexity class \emph{\textsf{NEXPTIME}}.
\end{mycorollary}

\section{Concluding remarks}\label{se:concluding remarks}
In this paper, we provided a positive solution to the s.p.\ for the slightly simplified variant \BSTuC of \MLSC, whose decision problem has been a long-standing open problem in computable set theory. \BSTuC differs from \MLSC in that membership has been dropped and the Cartesian product has been replaced by its unordered variant $\otimes$. Despite such simplifications, the s.p.\ for \BSTuC remains fully representative of the combinatorial difficulties due to the presence of the Cartesian product operator. Specifically, we proved that that both the ordinary s.p.\ and the (hereditarily) finite s.p.\ for \BSTuC are in \textsf{NEXPTIME}.

We conjecture that a more elaborated approach, inspired to the optimization results stated in Lemmas~\ref{was_Lemma5} and \ref{le:small fulfilling map} for the s.p.\ of \BST, will allow us to prove the \NP-completeness of the s.p.\ for \BSTuC.

We also expect that the technique introduce in this paper, based on $\otimes$-graphs and fulfilling maps, may be adapted to ascertain the decidability of various extensions of \BST with operators belonging to a specific class of operators, which includes, among others, the (ordered) Cartesian product $\times$ and the power set operator $\Pow$ and its variants $\powAst$ and $\powAstot$.

Finally, we are very confident that the decidability result for \BSTuC can be generalized to \MLSC, though at the cost of extra-technicalities, and we plan to report about it in a next paper.


\small

\begin{thebibliography}{99}

\bibitem[Can91]{Can91}
D.~Cantone.
\newblock Decision procedures for elementary sublanguages of set theory. {X}.
  {M}ultilevel syllogistic extended by the singleton and powerset operators.
\newblock {\em Journal of Automated Reasoning}, 7(2):193--230, 1991.

\bibitem[CCP90]{CCP90}
D.~Cantone, V.~Cutello, and A.~Policriti.
\newblock Set-theoretic reductions of {H}ilbert's tenth problem.
\newblock In {\em Proc. of 3rd Workshop ``Computer Science Logic'' 1989}, pages
  65--75, 1990.
\newblock Lecture Notes in Computer Science, 440.

\bibitem[CCS90]{CCS90}
D.~Cantone, V.~Cutello, and J. T. Schwartz.
\newblock Decision problems for Tarski's and Presburger's arithmetics extended with sets.
\newblock In {\em In E. B\"orger, H. B\"uning, M. Richter, and W. Sch\"onfeld, editors, Proceedings of 3rd Workshop Computer Science Logic - CSL ’90 (Heidelberg 1990)}, pages 95--109
  Berlin, 1990.
\newblock Lecture Notes in Computer Science, 533.

\bibitem[CDMO19]{CDMO19} D. Cantone, A. De Domenico, P. Maugeri, and E.G. Omodeo. Polynomial-time satisfiability tests for Boolean fragments of set theory. In A. Casagrande and E.G. Omodeo, editors, Proceedings of the 34th Italian Conference on Computational Logic, Trieste, Italy, June 19-21, 2019, volume 2396 of CEUR Workshop Proceedings, pages 123--137. CEUR-WS.org, 2019.

\bibitem[CFO89]{CFO89}
D.~Cantone, A.~Ferro, and E.G. Omodeo.
\newblock {\em Computable Set Theory}, vol.\ 6 
  {\em International Series of Monographs on Computer Science}.
\newblock Clarendon Press, Oxford, UK, 1989.

\bibitem[COP20]{COP20}
D.~Cantone, E.G.~Omodeo, and M.~Panettiere. From Hilbert's 10th problem to slim, undecidable fragments of set theory. In G. Cordasco, L. Gargano, and A. A. Rescigno, editors, \emph{Proceedings of the 21st Italian Conference on Theoretical Computer Science, ICTCS 2020, volume 2756 of CEUR Workshop Proceedings}, pages 47--60. CEUR-WS.org, 2020.

\bibitem[COP90]{COP90}
D.~Cantone, E.G.~Omodeo, and A.~Policriti. The automation of syllogistic. II: Optimization and complexity issues. {\em Journal of Automated Reasoning}, 6(2):173?187, 1990.

\bibitem[COP01]{COP01}
D.~Cantone, E.G. Omodeo, and A.~Policriti.
\newblock {\em Set {T}heory for {C}omputing - {\small From decision procedures
  to declarative programming with sets}}.
\newblock Monographs in Computer Science. {S}pringer-{V}erlag, New York, 2001.


\bibitem[COSU03]{COSU03}
D. Cantone, E.G. Omodeo, J.T. Schwartz, and P. Ursino.
\newblock Notes from the logbook of a proof-checker's project.
\newblock In N. Dershowitz, editor, {\em Verification: Theory and Practice
  (Essays Dedicated to Zohar Manna on the Occasion of His 64th Birthday)}, vol.
  2772 of {\em Lecture Notes in Computer Science}, pp.\ 182--207,
  Springer-Verlag, Berlin, 2003.

\bibitem[COU02]{COU02}
D. Cantone, E.G. Omodeo, and P. Ursino.
\newblock Formative processes with applications to the decision problem in set
  theory:~{I}.~{P}owerset and singleton operators.
\newblock {\em Information and Computation}, 172(2):165--201, 2002.

\bibitem[CU14]{CU14}
D. Cantone and P. Ursino.
\newblock Formative processes with applications to the decision problem in set
  theory: {II.} {P}owerset and singleton operators, finiteness predicate.
\newblock {\em Inf. Comput.}, 237: 215--242, 2014.

\bibitem[CU18]{CU18}
D. Cantone and P. Ursino.
\newblock An Introduction to the Technique of Formative Processes in Set
  Theory.
\newblock {\em Springer International Publishing}, 2018.

\bibitem[DPR61]{DPR61}
M. Davis, H. Putnam, and J. Robinson.
\newblock The decision problem for exponential Diophantine
equations.
\newblock {\em Annals of Mathematics}, 74(2): 425--436, 1961.

\bibitem[FOS80]{FOS80a}
A.~Ferro, E.G. Omodeo, and J.T. Schwartz.
\newblock {D}ecision procedures for elementary sublanguages of set theory. {I}:
  {M}ultilevel syllogistic and some extensions.
\newblock {\em Comm. Pure Appl. Math.}, 33:599--608, 1980.

\bibitem[Mat70]{Mat70}
Yu. V. Matiyasevich.  Enumerable sets are Diophantine (in Russian). \emph{Dokl. AN SSSR},
191(2), 278--282, 1970. Translated in: \emph{Soviet Math. Doklady}, 11(2), 354--358. Correction Ibid 11 (6), 1970, vi. Reprinted on pp. 269--273 in: \emph{Mathematical logic in the 20th century}, G. E. Sacks,
(Ed.), Singapore University Press and World Scientific Publishing Co., Singapore and River Edge, NJ, 2003.

\bibitem[Hil02]{Hilbert-02}
\newblock D. Hilbert. Mathematical Problems.
\newblock {\em Bulletin of the American Mathematical Society}, 8(10), 437--479, 1902.

\bibitem[OCPS06]{OCPS06}
E.G. Omodeo, D. Cantone, A. Policriti, and J.T. Schwartz.
\newblock A {C}omputerized {R}eferee.
\newblock In M. Schaerf and O. Stock, editors, {\em {R}easoning,
  {A}ction and {I}nteraction in {AI} {T}heories and {S}ystems -- {\small Essays
  dedicated to {L}uigia {C}arlucci {A}iello}}, vol. 4155 of {\em Lecture Notes
  in Artificial Intelligence}, pp.\ 117--139. Springer Berlin/Heidelberg, 2006.
  
\bibitem[OPT17]{OPT17}
E.G. Omodeo, A. Policriti, and A. Tomescu. On Sets and Graphs: Perspectives on Logic and Combinatorics. {\em Springer International Publishing}, 2017.

\bibitem[OS02]{OS02}
E.G. Omodeo and J.T. Schwartz.
\newblock A `{T}heory' mechanism for a proof-verifier based on first-order set
  theory.
\newblock In A.~Kakas and F.~Sadri, editors, {\em Computational Logic: Logic
  Programming and Beyond {\small-- Essays in honour of Bob Kowalski}, Part II},
  vol. 2048 of {\em Lecture Notes in Artificial Intelligence}, pp.\ 214--230.
  Springer-Verlag, Berlin, 2002.

\bibitem[Rob56]{Rob}
R.M. Robinson.
\newblock Arithmetical representation of recursively enumerable sets.
\newblock {\em Journal of Symbolic Logic}, 21(2), 162–186, 1956.

\bibitem[Sca84]{Sca}
B.~Scarpellini.
\newblock Complexity of subcases of Presburger Arithmetic.
\newblock {\em Transactions of the
American Mathematical Society}, 284(I):93--119, 1984.

\bibitem[Schw78]{Schw78}
J.T.~Schwartz.
\newblock Instantiation and decision procedures for certain classes of quantified set-theoretic formulae.
\newblock {\em ICASE Report}, 78-10, 1978.

\bibitem[SCO11]{SchCanOmo11}
J.T. Schwartz, D.~Cantone, and E.G. Omodeo.
\newblock {\em Computational logic and set theory: Applying formalized logic to
  analysis}.
\newblock Springer-Verlag, 2011.
\newblock Foreword by M. Davis.

\bibitem[Urs06]{Urs06}
P. Ursino.
\newblock A generalized small model property for languages which force the
  infinity.
\newblock {\em Matematiche (Catania)}, LX(I):93--119, 2005.


\end{thebibliography}
\end{document}